\documentclass[leqno,10pt]{article} 
\usepackage{amsmath, amssymb}
\usepackage{amsthm} 
\theoremstyle{plain} 

\theoremstyle{definition} 

\def\L{{\cal L}}

\def\bgn{\begin}

\def\L{{\cal L}}

\def\1{{[1]}}
\def\2{{[2]}}
\def\3{{[3]}}
\def\({\left(}
\def\){\right)}
\def\s-circ{\,{\scriptstyle{\circ}}\,}
\def\<<{<\negthinspace \negthinspace<}

\def\bgn{\begin}

\def\endaln{\end{align}}

\def\<{<\negthinspace \negthinspace <}

\def\({\left(}
\def\){\right)}

\def\[{\big[\neg\big[}
\def\]{\big]\neg\big]}
\def\al{\al}
\def\cpt{\text{\rm cpt}}
\def\M{{\cal M}}
\def\ind{\text{\rm ind}}
\def\diam{\text{\rm diam}}
\def\Vol{\text{\rm Vol}}
\def\tr{\text{\rm tr}}
\def\a{\alpha}
\def\b{\beta}

\def\e{\varepsilon}

\def\Gam{\Gamma}
\def\k{\kappa}
\def\del{\delta}
\def\lam{\lambda}
\def\ome{\omega}
\def\Ome{\Omega}
\def\sig{\sigma}

\def\R{\Bbb R}
\def\C{\Bbb C}

\def\O{\Bbb O}

\def\M{\frak M}

\def\w{\wedge}
\def\({\left(}
\def\){\right)}

\def\nab{\nabla}
\def\neg{\negthinspace}
\def\h{\hat}

\def\til{\tilde}

\def\ol{\overline}
\def\pa{\partial}
\def\olpa{\ol{\partial}}

\def\ran{\rangle} 
\def\lan{\langle}

\def\ss{\scriptscriptstyle}
\def\trian{\triangle}

\def\bsh{\backslash}

\def\:{\, :\,}

\def\K\"ahler{generalized K\"ahler}
\def\vol{\text{\rm vol}}
\def\10{\displaystyle L^{10}}
\def\2{\displaystyle L^2}
\def\c0{\displaystyle C^0}
\def\dstyle{\displaystyle}
\def\10{\displaystyle L^{10}}
\def\2{\displaystyle L^2}
\def\del{\delta}
\def\del2{\displaystyle L^2_{0,\delta}}
\def\c0{\displaystyle C^0}
\def\dstyle{\displaystyle}

\def\del{\delta}

\def\K{{\cal K}}
\def\M-A{\text{\rm Monge-Amp\`ere}}
\def\O{{\cal O}}
\def\cyl{\text{\rm cyl}}
\def\cone{\text{\rm cone}}
\def\Ric{\text{\rm Ric}}
\def\M-A{\text{\rm Monge-Amp\`ere}}

\makeatletter
\def\address#1#2{\begingroup
\noindent\parbox[t]{7.8cm}{%
\small{\scshape\ignorespaces#1}\par\vskip1ex
\noindent\small{\itshape E-mail address}%
\/: #2\par\vskip4ex}\hfill%
\endgroup}%
\makeatother
\title{{Calabi-Yau structures and Einstein-Sasakian structures \\
on crepant resolutions of isolated singularities
}} 
\author{
\textsc{Ryushi Goto$^{*}$} 
}
\date{} 
\begin{document}

\maketitle

\footnote{ 
$^{*}$Partly supported by the Grant-in-Aid for Scientific Research (C),
Japan Society for the Promotion of Science. 
}

\abstract{Let $X_0$ be an affine variety with only normal isolated singularity at $p$.
We assume that the complement $X_0\bsh\{p\}$ is biholomorphic to the cone 
$C(S)$ of an Einstein-Sasakian manifold $S$ of real dimension $2n-1$.
If there is a resolution of singularity $\pi : X\to X_0$ with trivial canonical line bundle $K_X$, then
there is a Ricci-flat complete K\"ahler metric for every K\"ahler class of $X$.
We also obtain a uniqueness theorem of Ricci-flat conical K\"ahler metrics in each K\"ahler class with a certain boundary condition. 
We show there are many examples of Ricci-flat complete K\"ahler manifolds arising as crepant resolutions.}

\tableofcontents
\section*{Introduction} 
Let $X$ be a K\"ahler manifold of complex dimension $n$ with trivial canonical line bundle $K_X$ and 
$\Ome$ a nowhere vanishing holomorphic $n$-form on $X$. 
If a K\"ahler form $\ome$ satisfies the following equation, 
$$
\Ome \w \ol{\Ome }=c_n\ome^n,
$$
for a constant $c_n$, 
then the Ricci curvature of $\ome$ vanishes, that is, $\ome$ is a Ricci-flat K\"ahler metric,
where $\ol\Ome$ is the complex conjugate 
of $\Ome$.
The well-known Calabi-Yau theorem, due to Yau on a compact K\"ahler manifold 
with the first Chern class $c_1=0$ was proved by solving the \M-A equation which shows that there exists a unique Ricci-flat K\"ahler metric in each K\"ahler class.
 On a non-compact complete K\"ahler manifold $X$ with $c_1(X)=0$, 
it is  an outstanding problem whether there exist Ricci-flat K\"ahler metrics. We need to impose suitable asymptotic conditions on the boundary.
There are many remarkable results on Calabi-Yau theorem on non-compact complete K\"ahler manifolds with $c_1=0$. 
Tian and Yau \cite{TY1}, \cite{TY2}, Bando-Kobayashi \cite{BK2} and Joyce \cite{Jo}  solved the \M-A 
equation under 
various boundary conditions. 
On the other hand, the hyperK\"ahler quotient construction 
\cite{HKLR}, \cite{Kr} 
produces many Ricci-flat K\"ahler manifolds
 in a simple and algebraic way, some of which are not obtained 
 by the analytic method \cite{Go1}.
Recently rapid developments occur in the Einstein-Sasakian geometry which yield a new view point of the problem of Ricci-flat K\"ahler metrics. 
If we have a positive Einstein-Sasakian manifold $S$, then the cone $C(S)=\Bbb R_{>0}\times S$  admits a Ricci-flat K\"aher cone metric. 
Boyer-Galicki \cite{BG} constructed a family of positive Einstein-Sasakian metrics on 
the links on hypersurfaces with isolated singularities, which includes interesting examples such as
homology spheres. 
Martelli and Sparks \cite{MS1} constructed explicitly a Sasaki-Einstein metric
on the sphere bundle of the canonical line bundle 
of the blown up $\C P^2$ at one point which is a irregular Sasaki manifold. 
Futaki, Ono and Wang \cite{FOW} showed that a sphere bundle of the canonical line bundle on every toric compact Fano manifold admits an Einstein-Sasakian metric. 
 These results imply that there are many notable examples of 
 Ricci-flat cone metrics which are constructed by Einstein-Sasakian manifolds.
 In the present paper, we introduce conical K\"ahler metrics which are complete K\"ahler metrics with the certain boundary condition 
 (see the definition\ref {def: conical metric} in section 1).
We apply an existence theorem of Ricci-flat K\"ahler metrics to the class of conical K\"ahler metrics. 
Let $\ome$ be a conical K\"ahler metric on $X$ with 
  $\Ome\w\ol\Ome = c_n F\ome^n$ for a positive function $F$. 
  If $F$ satisfies  $$\| e^{(2+\del)t}(F-1)\|_{C^k}<\infty,$$ for $0<\del<2n-2$ and $k\geq 2,$ then there exists a Ricci-flat conical K\"ahler metric $\ome_u$ on $X$ 
 (see the theorem \ref{th: existence CY on conical} in section 1 
for more detail). 
 
The existence theorem can be also deduced from the arguments by Bando-Kobayashi in \cite{BK2} and Tian-Yau \cite{TY2}. 
 Our boundary conditions in the theorem are modified for conical K\"ahler manifolds.
 The author gives a proof of the theorem for the sake of readers in section 2.
 
As an application, we discuss the existence of Ricci-flat K\"ahler metrics on 
resolution $X$ with trivial $K_X$ of an affine variety $X_0$ with only normal isolated singularity $\{p\}$.

\noindent
{\bf Theorem 5.1} \,\, {\it 
Let $X_0$ be an affine variety with only normal isolated singularity at $p$.
We assume that the complement $X_0\bsh\{p\}$ is biholomorphic to the cone 
$C(S)$ of an Einstein-Sasakian manifold $S$ of real dimension $2n-1$.
If there is a resolution of singularity $\pi : X\to X_0$ with trivial canonical line bundle $K_X$, then
there is a Ricci-flat complete K\"ahler metric for every K\"ahler class of $X$. }

\par\medskip

Our theorem covers the crucial case where a K\"ahler class does not belongs to the compactly supported cohomology group. 
We use a vanishing theorem on $X$ and the Hodge and the Lefschetz decomposition theorems on a Sasakian manifold to construct a suitable initial K\"ahler metric in every K\"ahler class which the existence theorem can be 
applied.
We show that a Ricci-flat K\"ahler conical metric is unique in
 each K\"ahler class if we impose a certain boundary condition on metrics.
(see the theorem \ref{th: uniqueness in class}).\par
In section $1$, we introduce the class of conical K\"ahler metrics and 
show the existence theorem and the uniqueness theorem of Ricci-flat K\"ahler metrics on them. 
In section 2  we give a proof of the existence theorem and the uniqueness theorem.
In section 3, we will give a short explanation of Sasakian metrics
and K\"ahler cone metrics. 
In section 4 we discuss a one to one correspondence between Einstein-Sasakian structures and Ricci-flat K\"ahler cone metrics. 
 In section $5$, as an application of section $1$, 
 we obtain Ricci-flat K\"ahler conical metrics on crepant resolutions of 
 normal isolated singularities as above. 
 In section $6$, we construct several families of Ricci-flat K\"ahler conical metrics on crepant resolutions of 
 normal isolated singularities.
 Our examples include: 
  resolutions of the isolated quotient singularities,
 the total spaces of canonical line bundles of K\"ahler-Einstein Fano manifolds, 
 the total space of the canonical line bundle of every toric Fano manifold and 
 small resolutions of ordinary double points of dimension $3$ 
  \par 
Some of these examples are already known. Joyce \cite{Jo} showed the Calabi-Yau theorem on resolutions of the isolated quotient singularities which is called asymptotically locally Euclidean (ALE). 
 Calabi \cite{Cal} used the bundle construction to obtain Ricci-flat K\"ahler metrics on 
 the total spaces of the canonical line bundles of K\"ahler-Einstein Fano manifolds 
It must be noted that the K\"ahler classes of these Ricci-flat K\"ahler metrics
 lie in the compactly support cohomology group.
 Van Coevering \cite{Coe} and Santoro \cite{San} constructed Ricci-flat K\"ahler metrics
on crepant resolutions whose K\"ahler classes also belong to the compactly support cohomology group.

Our method provides wider classes of complete Ricci-flat K\"ahler metrics
(see the examples in section 6)  and the theorem 5.1 shows that the conjecture discussed in \cite{MS} and \cite{Coe} on 
the existence of complete Ricci-flat K\"ahler metrics on resolutions of cones is affirmative.

After the author posted this paper in the Arxiv, Van Coevering \cite{Coe2} submitted a paper 
which covers the existence theorem 1.5 by the different method. 

The author would like to thank Professor A. Futaki and R. Kobayashi 
for suggestive discussions about Ricci-flat K\"ahler metrics. 
He thanks Professor S. Bando for his kind and remarkable comments about 
the existence theorem in \cite{BK2}.
He is grateful to Prof. A. Fujiki for his helpful advises. 

\numberwithin{equation}{section}
\section{Existence theorem of Ricci-flat conical K\"ahler metrics}
Let $(S, g_S)$ be a compact Riemannian manifold of dimension $2n-1$ and 
$C(S)=\R_{>0}\times S$ the product of the positive real number $\R_{>0}$ and $S$
with a coordinate $r\in \R_{>0}$, which is called the cone of $S$.
By changing the coordinate  $t=\log r$, 
we regard the cone $C(S)$ as the cylinder $\R\times S$ with the cylinder parameter $t\in \R$ 

\bgn{definition}\label{def: cyl and cone metric}
The cylinder metric $g_{\cyl}$ on $C(S)$ is the product metric,
$$dt^2+ g_S$$
and the cone metric $g_{\cone}$
is given by 
$$
g_{\cone}= dr^2+ r^2 g_S.
$$
Let $\nab_{cyl}$ be the Levi-civita connection with respect to the cylinder metric $g_{\cyl}$ and 
$| \a|_{g_{\cyl}}$ the point-wise norm of a tensor $\a$ by $g_{\cyl}$. 
The $C^k$-norm of a tensor $\a$ is given by 
$$
\|\a\|_{C^k} =\sum_{i=0}^k \sup |\nab_{\cyl}^i \a|_{g_{\cyl}}.
$$
We also have the H\"older norm 
$\|\a\|_{C^{k,\a}}$, for $0<\a<1$.
In this paper we use the $C^k$-norm and 
$C^{k,\a}$-norm with respect to the cylinder metric $g_{\cyl}$ unless it is  mentioned. 
\end{definition}
Since $r=e^t$, we have the relation
\bgn{equation}
g_{\cone}=r^2g_{\cyl}
\end{equation}
which follows from 
\bgn{align*}
r^2g_{\cyl}=&r^2 (dt)^2+r^2g_S\\
=&(dr)^2+r^2 g_S \\
=&g_{\cone}.
\end{align*}
\bgn{definition} A manifold $X$ has {\it a cylindrical boundary} if there is a compact set $K$ of $X$ such that the complement $X\bsh K$ is diffeomorphic to the cylinder $C(S)=\R\times S$. 
We identify $X\bsh K$ with the cone $C(S)$.
A Riemannian metric $\ol{g}$ on a manifold $X$ with a cylindrical boundary $C(S)$ is {\it a cylindrical metric} 
if $\ol{g}$ satisfies the following condition on $X\bsh K\cong C(S)$,
$$
\| e^{\del t}(\ol{g}-  g_{\cyl}) \|_{C^k}<\infty, 
$$
for some $\del>0$ and an integer $k>4$.
In other words, the difference between  $\ol{g}$ and $g_{\cyl}$ decays exponentially with order $O(e^{-\del t})$, including their higher order derivatives up to $k$
$$
\sum_{i=0}^k |\nab_{\cyl}^i(\ol{g}- g_{\cyl})|_{g_{\cyl}}= O(e^{-\del t}),
$$
\end{definition}
We always extend the cylinder parameter $t$ as a $C^\infty$ 
function on $X$ and identify the complement $X\bsh K$ with the cone $C(S)$.
\bgn{definition}\label{def: conical metric}
A Riemannian metric $g$ on a manifold with a cylindrical boundary $C(S)$ is {\it a conical metric} if $r^{-2}g=e^{-2 t}g$ is
a cylindrical metric, that is, $g$ satisfies 
$$
\| e^{(-2+\del) t}(g-  g_{\cyl}) \|_{C^k}<\infty, 
$$
A conical K\"ahler metric on a complex manifold $X$ is a conical Riemannian metric which is K\"ahlerian.
A conical K\"ahler form $\ome$ is a K\"ahler form with the associate Riemannian metric is conical and 
a manifold with a conical K\"ahler form is called a conical K\"ahler manifold.
\end{definition}
\bgn{definition}
Let $(X,\ome)$ be a conical K\"ahler manifold with trivial canonical line bundle $K_X$ 
and $\Ome$ a nowhere vanishing holomorphic $n$-form $\Ome$ on $X$. 
If a pair $(\Ome,\ome)$ satisfies the following equation 

$$
\Ome\w\ol\Ome = c_n \ome^n,
$$
for a constant $c_n$, 
then $(\Ome, \ome)$ is called a Calabi-Yau structure
whose K\"ahler form $\ome$ 
gives the Ricci-flat K\"ahler metric.
\end{definition}

The following theorem can also deduced from the arguments by Bando-Kobayashi \cite{BK2} and Tian-Yau \cite{TY2}.
For the completeness of the paper, the author gives a proof of existence theorem in section 2. 
\bgn{theorem}\label{th: existence CY on conical}
Let $(X,\ome)$ be a conical K\"ahler manifold of complex dimension $n$ with trivial canonical line bundle $K_X$ 
and $\Ome$ a nowhere vanishing holomorphic $n$-form $\Ome$ on $X$ which defines a positive function 
$F$ by 
$$
\Ome\w\ol\Ome = c_n F\ome^n. 
$$
If $F$ satisfies the following condition 
$$
\| e^{(2+\del)t} (F-1)\|_{C^{k,\a}} <\infty, 
$$
for $0<\del < 2n-2$ and $k\geq 2$, $0<\a<1$,
then there exists a smooth solution $u$ of the \M-A equation
$\Ome\w\ol\Ome =c_n \ome_u^n$, such  that
\bgn{equation}\label{eq: ome u}
\ome_u= \ome +\sqrt{-1}\pa\ol\pa u 
\end{equation} is a conical Ricci-flat K\"ahler form with the condition,
$$
\| e^{\del t} u\|_{C^{k+2,\a}} <\infty.
$$

(Note that $\ome_u$ is a complete Ricci-flat K\"ahler metric.)
\end{theorem}
\bgn{remark}
In the theorem \ref {th: existence CY on conical}, the decay order of $F-1$ is crucial for the existence 
of Ricci-flat K\"ahler metrics.
Note that we estimate it by the cylinder metric $g_{\cyl}$.
In order to obtain a solution of the \M-A equation, 
we need to solve the equation of the Laplacian with respect to the conical metric $g$,
$\trian_g u =v$ on $X$, which is a linearization of the \M-A equation.
The equation $\trian_g u =v$ has a unique solution $u$ with the decay order 
$O(e^{-\del t})$ if $v$ decays with order $O(e^{-(2+\del)t})$
for $0<\del<2n-2$.  
Thus we require that $F$ decays with  the order $O(e^{-(2+\del) t})$. 
\end{remark}

We show that a Ricci-flat conical K\"ahler metric in the form as in 
theorem \ref {th: existence CY on conical} is unique.
\bgn{theorem}\label{th: uniqueness theorem }
If there are two Ricci-flat conical K\"ahler metrics $\ome$ and $\ome'$ satisfying 
$$
\ome'=\ome+\sqrt{-1}\pa\ol\pa u,
$$
where $u$ is a function with $\|e^{\del t}u\|_{C^k}<\infty$, for a positive $\del$ and $k\geq 2$.
Then $\ome=\ome'$. 
\end{theorem}
There is the action of the automorphism group of $X$ on Ricci-flat K\"ahler metrics, in this sense, the uniqueness theorem does not hold on 
non-compact K\"ahler manifolds. 
However if we impose the following boundary condition on metrics, we obtain the uniqueness theorem of Ricci-flat conical K\"ahler metrics in each K\"ahler class 
\bgn{theorem}\label{th: uniqueness in class}

Let $\ome$ and $\ome'$ be two Ricci-flat conical K\"ahler metrics on $X$ of dimension $n$ 
with $[\ome]=[\ome']\in H^2(X)$ ($n\geq 2$).
We assume that $X$ has a cylindrical boundary $C(S)=\R\times S$ with 
$H^1(S)=\{0\}$.
If  $\| e^{\lam t}(\ome-\ome')\|_{C^k}<\infty$, 
for a constant $\lam >n-2$ and $k\geq 2$, then $\ome=\ome'$.
\end{theorem}

\section{Proof of the existence theorem}
\subsection{The Laplacian on conical Riemannian manifolds}
Let $(X, \ol{g})$ be a $2n$ dimensional cylindrical Riemannian manifold with cylindrical boundary 
$C(S)$ and $g_{\cyl}$ the cylinder metric on $C(S)=\R\times S$ as in section $1$, 
$g_{\cyl} =dt^2+ g_S$, where $t\in \R$ is the cylinder parameter on $C(S)$ and $g_S$ is a Riemannian metric on the manifold $S$ of dimension $2n-1$.
We denote by $\| f\|_{L^p_k}$ the Sobolev norm of a function $f$ on $X$ with respect to the cylindrical metric $\ol{g}$, 
$$
\|f\|_{L^p_k}=\sum_{i=0}^k \(\int_X |\nab^i f|^p\vol_{\cyl}\)^{\frac1p},
$$
where $\vol_{\cyl}$ is the volume form with respect to $\ol{g}$.
We define a weighted Sobolev norm on $X$ by using the exponential function
$e^{\del t}$, 
\bgn{equation}
\|f\|_{L^p_{k,\del}} =\| e^{\del t} f\|_{L^p_k}
\end{equation}
The weighted Sobolev space is the completion of $C^\infty$ functions on $X$ with compact support with respect to the weighted Sobolev norm.
Note that a cylindrical Riemannian manifold is complete.
We also define a weighted H\"older norm as 
$\| f\|_{C^{k,\a}_\del}= \|e^{\del t} f\|_{C^{k,\a}}$ for 
$k$ and $0<\a<1$ with respect to the cylindrical metric.
Recall that a Riemannian metric $g$ is conical if $r^{-2} g= e^{-2 t} g$ is a cylindrical metric on $X$. 
On a conical Riemannian manifold $(X, g)$, we use the weighted Sobolev norm 
$\|f\|_{L^p_{k,\del}}$ with respect to the cylindrical Hermitian metric $r^{-2} g$ and 
the weighted H\"older norm is also defined in terms of the cylindrical metric $r^{-2} g$.
If a function $f$ lies in $C^{k,\a}_\del$ implies that $f$ decays exponentially 
with order $O(e^{-\del t})$ together with its derivatives. 
(We call $\del$ the weight.) Then the Laplacian $\trian_g$ gives the bounded linear operator from $L^p_{k+2,\del}$ to $L^p_{k,\del+2}$.
\bgn{remark}
Let $\vol_{\cone}$ be the volume form with respect to the conical metric $g$.
Then note that
$$
\vol_{\cone}= r^{-2n}\vol_{\cyl}.
$$
\end{remark}
Let $*_g$ be the Hodge star operator with  respect to the conical metric $g$ and 
$*_{\cyl}$ the Hodge star operator with respect to the cylindrical metric $r^{-2} g$.
\bgn{lemma}
We define a differential operator $P$ by 
$$P=r^2\trian_g.$$  Then the operator $P$ is given by 
\bgn{equation}\label{eq: r 2triang }
P=r^2\trian_g=\trian_{\cyl}-(2n-2)\frac{\pa}{\pa t}
\end{equation}
on the cylinder boundary $C(S)=X\bsh K\cong \R\times S$, 
where $t$ is the cylinder parameter and $\trian_{\cyl}$
 the Laplacian with respect to the cylindrical metric $r^{-2} g$.
\end{lemma}
\bgn{proof}
 Since $\trian_g= *_gd*_gd$ and $\trian_{\cyl}=*_{\cyl} d{*_{\cyl}}d$,
comparing two Hodge star operators 
 $*_g$ and $*_{\cyl}$, we have 
(\ref{eq: r 2triang }).
\end{proof}\par
The cylinder metric $g_{\cyl}=dt^2+g_S$ on $C(S)$ gives the Laplacian 
$\trian_{g_{\cyl}}=(-i\frac{\pa}{\pa t})^2+\trian_S$ which is invariant under the translation of cylinder parameter $t$, where 
$\trian_S$ is the Laplacian of $(S, g_S)$.

Let $I(P)$ be the translation-invariant operator given by 
\bgn{align}
I(P) =& \trian_{g_{\cyl}}-(2n-2)\frac{\pa}{\pa t}\\
=&(-i\frac{\pa}{\pa t})^2+\trian_S-(2n-2)\frac{\pa}{\pa t}.\label{eq: I(P)}
\end{align}
Then the all coefficients of the operators $P -I(P)$ decays exponentially together with their 
derivatives, 
which implies that we can apply the theory of elliptic differential operators on cylindrical Riemannian manifolds developed by \cite{LMc},  \cite{Mel}. 
(Note that these operators are called $b$-operators in \cite{Mel}.)
We can select suitable weighted Sobolev spaces for the operator $P$ to be a Fredholm operator.
Substituting $\lam$ for $-i\frac{\pa}{\pa t}$ into (\ref{eq: I(P)}),
we define a family of the operators $I(P,\lam)$ on $S$ given by  
$$
I(P,\lam ) =\lam^2-(2n-2)\lam i +\trian_S
$$
which is called the indicial family parametrized by $\lam \in \C$.
The operator 
$$
I(P,\lam): L_{k+2}^p(S)\to L_k^p(S)
$$
is an isomorphism for all $\lam\in  \C\bsh $Spec$(P)$, that is, 
The operator $I(P)$ does not admit a bounded inverse for $\lam\in $Spec$(P)$.
The properties of the set Spec$(P)$  are  developed in \cite{LMc} and \cite{Mel}. 
In our cases, 
Spec$(P)$ is explicitly described as 
\bgn{lemma}
$$
\text{\rm Spec}(P) = \{ \,0, \,\,\,(2n-2)\sqrt{-1} , \,\,\,\,\mu_j^+\sqrt{-1},\,\, \,\,\mu_j^-\sqrt{-1} \,\,|\, j=1,2,\cdots\}
$$
where
$\mu_j^\pm\sqrt{-1}$ are two solutions of the quadratic equation
  $x^2-(2n-2)\sqrt{-1}\,x+\lam_j=0$
and $\lam_j$ is the $j$-th eigenvalue of the Laplacian $\trian_S$ on $S$.
\end{lemma}
\bgn{proof}
It follows from $I(P,\lam ) =\lam^2-(2n-2)\lam i +\trian_S $.
\end{proof}\par
Let ImSpec$(P)$ be the set consisting of the imaginary parts of elements of Spec$(P)$. 
Then Imspec$(P)$ is given by 
$$
\cdots < \mu_2^-< \mu_1^-< 0< 2n-2 <\mu_1^+< \mu_2^+
$$
We select the weighted Sobolev space
$L^p_{k,\del}$ with weight $\del\notin$ImSpec$(P)$. 
Then from theorem 1.1 in \cite{LMc} and \cite{Mel},
$P: L^p_{k+2,\del}\to L^p_{k,\del}$ is a Fredholm operator with index 
$\ind(P, \del)$. 
Further for  $0<\del<n-1$, since the Laplacian $\trian_g$ is self-adjoint, applying the theorem 
7.4 (see page 436) in \cite{LMc}, we obtain 
$$
\ind(P, \del+n-1) +\ind(P,-\del+n-1)=0
$$
There is the formula for the index of the  cylindrical operators  
depending on weights $\del$ in \cite{LMc}.
Applying the theorem 1.2 in \cite{LMc}
(Note that $N(-\del+n-1,\del+n-1)=0$), we have 
$$
-\ind(P, \del+n-1) +\ind(P,-\del+n-1)=0.
$$
Thus we obtain
$\ind(P, \del)=0$ for  $0<\del<2(n-1)$.

By the result, we have the solvability of the equation of the Laplacian
 $\trian_g$.
\bgn{proposition}\label{prop: Laplacian(conical) Lp}
Let $(X, g)$ be a $2n$-dimensional conical Riemannian manifold 
 and $\trian_g$ the Laplacian with respect to the conical metric $g$.
For a weight $\del$ with $0< \del <2n-2$, 
the Laplacian $\trian_g$ gives an isomorphism between 
$L^p_{k+2,\del}$ and $L^p_{k,\del +2}$, for positive integers $p$, $k$ with 
$\frac 1p< \frac{k-\a}{2n}$ for $0<\a<1$.
\end{proposition}
\bgn{proof}
It follows from the Sobolev embedding theorem that 
$L^p_{k+2,\del}\subset C^{2,\a}_\del$. 
By the Maximum principle, a harmonic function attains its Maximum and Minimum at the boundary. 
Since a function in $C^{2,\a}_\del$ decays exponentially, we have 
$\ker\trian_g=\{0\}$ which implies that 
$\ker P =\ker \trian_g=\{0\}$.
Since $\ind(P, \del) =0$, we have an isomorphism 
$P: L^p_{k+2, \del}\to L^p_{k,\del}$. Since $\trian_g$ is the composition $r^{-2}P$,
$\trian_g: L^p_{k+2,\del}\to L^p_{k, \del+2}$ is an isomorphism also.
\end{proof}\par
Then we have the followings,
\bgn{proposition}\label{prop: Laplacian(conical)}
Let $(X, g)$ be a $2n$-dimensional conical Riemannian manifold 
 and $\trian_g$ the Laplacian with respect to the conical metric $g$.
For a weight $\del$ with $0< \del <2n-2$, 
the Laplacian $\trian_g$ gives an isomorphism between 
$C^{k+2,\a}_{\del}$ and $C^{k,\a}_{\del +2}$.
In other words, there is a unique solution $u \in C^{k+2, \a}_\del$ of the equation 
$\trian_g u =f$ for all $f \in C^{k,\a}_{\del+2}$, where $k\geq 0$.
\end{proposition}
\bgn{remark}
If $(X,g)$ is an ALE space, then the proposition
\ref{prop: Laplacian(conical)} implies the theorem 8.3.5 (a)
in \cite{Jo}. 
The H\"older space $C^{k,\a}_\b$ in \cite{Jo} coincides with 
our H\"older space $C^{k,\a}_{\del}$ with $\del =-\b$.
\end{remark}
We need the following lemma for a proof of the proposition \ref{prop: Laplacian(conical)}.
\bgn{lemma}\label{lem: maximum principle}
Let $(X, g)$ be a conical Riemannian manifold with the cylinder parameter $t$ and 
$\trian_g$ the Laplacian with respect to the conical metric $g$. 
We take a smaller weight $\til\del$ with $0<\til\del<\del<2n-2$.
We assume that a function $u\in C^k_{\til \del}$ satisfies 
$$
\trian_g u =h,
$$ 
and we already have a bound of $C^0$-norm $\|u\|_{C^0}<C_1$ and 
a weighted $C^0$-norm 
$\| e^{(2+\del) t} h\|_{C^0 }<C_2$.
Then there is a constant $C>0$ depending only on $C_1$, $C_2$ and $g$ such that 
$$
\|e^{\del t}u(x)\|_{C^0}< C
$$
\end{lemma}
\bgn{proof}{\it of lemma }\ref{lem: maximum principle}\,
Let $\trian_{\cone}$ be the Laplacian with respect to the cone metric $g_{\cone}$ on the cone $C(S)$. 
Then we see that 
$$
\trian_{\cone} = e^{-2t}\( -(\frac{\pa}{\pa t})^2 -( 2n-2) \frac{\pa}{\pa t} +\trian_S\)
$$
Thus we have 
$$
\trian_{\cone}e^{-\del t} =\del(2n-2-\del) e^{-(2+\del)t}>0
$$
on the cone $C(S)$, for $0<\del<2n-2$.
Then from the definition of the conical metric, there are constants $C_3>0$ and $T>0$ such that 
\bgn{equation}\label{eq: Maximum principle-1}
\trian_ge^{-\del t}\geq C_3e^{-(2+\del)t},
\end{equation}
on the region $D_T:=\{ t> T\}$.
We take a constant $C$ satisfying the two inequalities,
\bgn{align}
&C_3 C> C_2\label{eq: Maximum principle-2}\\
&Ce^{-\del T} -C_1>0.\label{eq: Maximum principle-3}
\end{align}
Then from (\ref{eq: Maximum principle-1}) and (\ref{eq: Maximum principle-2}), we have 
\bgn{align}
\trian_g( Ce^{-\del t}\pm u) >& CC_3e^{-(2+\del)t}\pm h\\
>&CC_3e^{-(2+\del )t}-C_2e^{-(2+\del)t}\\
>&0\label{eq: Maximum principle-4}
\end{align}
From (\ref{eq: Maximum principle-3}), we also have an inequality on the compact set $\{t\leq T\}$
\bgn{equation}\label{eq: Maximum principle-5} 
(Ce^{-\del t}\pm u) (x) \geq Ce^{-\del T} -C_1\geq 0, 
\end{equation}
for all $ x \in \{ t\leq T\}$.
If the function $Ce^{-\del t}\pm u$ have its Minimum at $x_0\in D_T$, 
then $\trian_g (Ce^{-\del t}\pm u)(x_0)\leq 0$. 
Thus from (\ref{eq: Maximum principle-4}), $Ce^{-\del t}\pm u$ 
can not have the Minimum on $D_T$.
Since $\lim_{t\to\infty} u(x)=0$, 
It follows from (\ref{eq: Maximum principle-5}) that 
$(Ce^{-\del t}\pm u)(x)\geq 0$, for all $x\in X$.
Hence we have $\|e^{\del t} u\|_{C^0}<C$.
\end{proof}\par

\bgn{proof}{\it of proposition }\ref{prop: Laplacian(conical)}\,
Since $u\in C^{k+2,\a}_\del$ decays exponentially, it follows that $\ker \trian_g=\{0\}$. 
For a smaller weight $0<\til\del<\del$, any function $f\in C^{k,\a}_{2+\del}$ is a function in $L^p_{k,\til\del+2}$. 
Then it follows from the proposition \ref{prop: Laplacian(conical) Lp} that there is a function $u \in L^p_{k+2,\til\del}$ such that 
$\trian_gu =f$.
Since $u$ is a function in $C^{k+2,\a}_{\til\del}$, it follows from the lemma \ref{lem: maximum principle} that 
$u\in C^0_\del$. 
We have the following equation, 
$$
\trian_g (e^{\del t} u) = H,
$$
where $H\in C^{k,\a}_\del$ is a function consisting of the terms $(\trian_ge^{\del t})u$, $e^{\del t}f$ and $(\nabla_g e^{\del t})(\nabla_g u)$.
Thus from the Schauder estimate, we have  
$e^{\del t} u \in C^{k+2,\a}$. 
Hence we have a unique solution $u\in C^{k+2, \a}_\del$ of the equation $\trian_gu =f$, for $f\in C^{k+2,\a}_\del$.
\end{proof}\par


\subsection{The Sobolev inequality for conical Riemannian manifolds}
\bgn{proposition}\label{prop: Sobolev (conical)}
Let $(X, g)$ be a conical Riemannian manifold of dimension $2n$.  
Then every function $f \in L^{2\e}\cap L^2_1$ satisfies: 
$$
\|f\|_{\displaystyle L^{2\e}} <S \| df\|_{\dstyle L^2},
$$
where $\e=\frac  n{n-1}$ and $S$ is a constant which depends only on $n$ and $g$.
Note that the the point-wise norm and $L^2$-norm are defined by the cone metric $g$ in this proposition. 
\end{proposition}
Note that the curvature of a conical Riemannian manifold is bounded. 
 In the proposition, we do not assume the condition: 
$$\int_X f\vol=0,$$ 
which is necessary for the Sobolev inequality for compact Riemannian manifolds.
 The theorem is already known for geometers, 
 see propositoin 3.2 in \cite{TY2}, however for the sake of readers, the author gives a proof. 
Our proof relies on a well-known result on compact Riemannian manifolds :
\bgn{theorem}\label{th: Sobolev constant}\text{\rm \cite{Gal}}
Let $M$ be a $2n$ dimensional compact Riemannian manifold on which the Ricci curvature $\Ric(g)$ satisfies 

$$
\text{\rm diam}(X,g)^2\text{\rm Ric}(g)\geq -\a g
$$
for a constant $\a$.
Then there is a constant $\k(2n,\a)$ depending only on 
 $2n$ and constant $\a$ such that for all $C^\infty$ function $f$ on $X$
 with $\int_X f\vol=0$,
the following Sobolev inequality holds,
$$
\k(2n,\a)\frac{\text{\rm Vol}(M,g)^\frac 1{2n}}{\text{\rm diam}(X,g)}\|f\|_{L^{2\e}}
\leq \|df\|_{L^2}
$$ 
\end{theorem}
\bgn{proof}(Proposition \ref{prop: Sobolev (conical)})
Since $(X,g)$ is complete, the Hilbert space $L^2_1\cap L^{2\e}$ is the completion of the space of 
compactly supported $C^\infty$ functions.
Hence it suffices to show that the theorem holds for compactly supported functions.
We extend the cylinder parameter $t$ as a $C^\infty$ function on $X$ and put $r=e^t$.
We define a compact subset $X_R$ by 
$$
X_R=\{ x\in X\, |\, r(x) \leq R \},
$$  
for a constant $R$.
We take two copies $X_R^1$, $X_R^2$ of $X_R$, 
with the cylinder parameter $t_i$ on each $X_R^i$ for $i=1,2$ and put $r_i=e^{t_i}$.
we obtain a compact manifold $M_R$ to glue 
$X_R^1$ and $X_R^2$  at the region $\{ x\in X^i_R\, |\, R-1\leq r_i(x)\leq R\}$ by the following diffeomrphism 
$\psi $, 
$$
\psi(r_1, y) = (-r_1+2R-1, y),
$$
where $y\in S$.
The metirc $g$ is conical and we have 
$$
g_i=dr_i^2+r_i^2g_S + a_i,
$$
where $a_i$ decays exponentially including their higher derivatives, i.e., $O(e^{-\del t})$,
for $i=1,2$. 
Then we have 
\bgn{align}
|\psi^*g_2-g_1|\leq&C_1|dr_1^2+(-r_1+2R-1)^2g_S - (dr_1^2+r_1^2g_S)|\\
\leq &C_1|(-2r_1(2R-1)+(2R-1)^2)g_S|\leq C_1|(2R-1)g_S |,
\end{align}
where $R-1\leq r_1\leq R$ and we take $R$ sufficiently large and $C_1$ is a constant which does not depend on $R$. 
(In this poof, we use $C_i$ as a constant which does not depend on $R$.)
Thus using the $C^2$ norm with respect to the cone metric $g_1$, we obtain 
\bgn{equation}\label{eq: proof of Sobolev inequality 1}
|\psi^*g_2- g_1|_{C^2}\leq C_2 r_1^{-1} 
\end{equation}
Note that $|g_S|=O(r^{-2})$ by the cone metric.
We define a smooth metric $g_R$ on $M_R$ by 
$g_R=\psi^*g_2+\chi_R(-\psi^*g_2+g_1)$
where $\chi_R$ is the cut-off function which is a smooth function with $0\leq \chi_R\leq 1$ and 
$$\chi_R=
\bgn{cases}
&1, \, ( r_1<R-1 ) \\
&0, \,(r_1>R)
\end{cases}
$$
Then $g_R$ is 
$$g_R=
\bgn{cases}
&g_1,\, ( r_1<R-1)\\
&g_2, \,(r_1>R)
\end{cases}
$$
Then It follows from (\ref{eq: proof of Sobolev inequality 1}) that 
$|g_R-g_2|_{C^2}\leq C_3r_1^{-1}$.
Thus we have a bound of  the curvature $R_{g_R}$ of $g_R$ by the curvature $R_{g_2}$ of $g_2$, 
$$
C_4^{-1}R_{g_2}< R_{g_R}< C_4R_{g_2}
$$
Since $g_2$ is conical metric, we have 
diam $(M_R, g_R)^2$ Ric$_{g_R}> -C_5 g_R$.
Since $X$ is a conical Riemannian manifold, the volume $\vol(M_R)$ and the diameter diam$(M_R)$ of 
$M_R$ satisfy 
$$
\frac{\diam(M_R)}{\Vol(M_R)^\frac 1{2n}}<C_5
$$
for a constant $C_5$ which does not depend on $R$.
For a compactly supported $C^\infty$ function $f$, we take $R$ sufficiently large such that 
the support of $f$ is a subset of $X_{R-1}$.
We define a function $f_R$ on $M_R$ by 
$$
f_R=
\bgn{cases}
&+f(x), \quad x\in M_R^1 \\
&-f(x),\quad  x\in M_R^2
\end{cases}
$$
Since 
$\int_{M_R} f_R \vol =0$, we apply the theorem \ref{th: Sobolev constant} to 
the function $f_R$ and we have 
$$
\| f_R\|_{L^{2\e}}<\k(2n,a)^{-1}C_5\|df_R\|_{L^2}
$$
Since we have 
\bgn{align}
&\|f_R\|_{L^{2\e}}=2\|f\|_{L^{2\e}}, \\ 
&\|df_R\|_{L^2}= 2\|df\|_{L^2},
\end{align}
and $C_5$ does not depend on $R$,
we have the inequality
$$\|f\|_{L^{2\e}}\leq \k(2n,\a)^{-1}C_5\|df\|_{L^2},
$$
for a $C^\infty$ function $f$ with compact support.
Thus we obtain the result.
\end{proof}

\subsection{The inequality of solutions of the \M-A equation}
Let $(X, g, \ome)$ be a conical K\"ahler manifold. 
We use the same notation as in the theorem \ref{th: existence CY on conical}.
We assume that 
there is a $C^\infty$ function $u$ on $X$ which satisfies the equation
$$
\ome_u^n=F\ome^n,
$$
where 
$$
\ome_u=\ome+dd^cu
$$
is a K\"ahler form on $X$ and 
a $C^\infty$ positive function $F$  satisfies 
$F-1 \in C^{k,\a}_{2+\del}$ and 
$u\in C^{k+2,\a}_\del$ and $d^c=\frac1{2\sqrt{-1}}(\pa-\ol\pa)$.
We select a weight  $\del$  with $0<\del<2n-2$
and a natural number $p>2$  with
$\del p +2>2n$.
Then since $\ome^n= e^{2nt}\vol_{\cyl}$,
the function $u|u|^{p-2}(1-F)$  is integrable with respect to the volume form $\ome^n$.
Thus we have 
\bgn{align*}
\int_X u|u|^{p-2}(1-F)\ome^n=&
\int_Xu|u|^{p-2}(\ome^n-\ome_u^n)\\
=\int_X u|u|^{p-2}(-dd^cu)\w
(\ome^{n-1}&+\cdots+\ome^{n-1}_u)
\end{align*}
Since $\del p +2>2n$ and 
$ u|u|^{p-2}d^c u\w(\ome^{n-1}+\cdots+\ome^{n-1}_u)$
decays exponentially with respect to the cylinder parameter $t$, 
we can apply the Stokes theorem,
$$
\int d\( u|u|^{p-2}d^c u\w(\ome^{n-1}+\cdots+\ome^{n-1}_u)\)=0
$$
Substituting 
$d(u|u|^{p-2})=(p-1)|u|^{p-2}du$, we obtain 
\bgn{align}\label{60}
\int_X u|u|^{p-2}(1-F)\ome^n=&(p-1)\int_X|u|^{p-2}du\w d^cu\w
(\ome^{n-1}+\cdots+\ome^{n-1}_u)
\end{align}
Since 
$du\w d^c u\w (\ome^i\w\ome_u^{n-1-i})\geq 0$ 
at each point on $X$,
we have 
$$
du\w d^c u\w(\ome^{n-1}+\cdots +\ome_u^{n-1})
\geq du\w d^cu\w\ome^n
$$
Substituting it into (\ref{60}), we have 
\bgn{align}
\int_X u|u|^{p-2}(1-F)\ome^n\geq&(p-1)\int_X|u|^{p-2}du\w d^cu\w\ome^{n-1}\\
=&\frac{p-1}n\int_X|u|^{p-2}|du|^2_g\ome^n\\
=&\frac{4(p-1)}{p^2n}\int_X  \big|(d\,|u|^{\frac p2})\big|_g^2\ome^n
\end{align}
where 
we are using  
$$
\frac14 p^2|u|^{p-2}|du|^2_g=\big|(d\,|u|^{\frac p2})\big|_g^2
$$
and  $\big|(d\,|u|^{\frac p2})\big|_g^2$ is the norm of the $1$-form $d\,|u|^{\frac p2}$ with respect to the conical metric $g$.  Hence we obtain the inequality,
\bgn{proposition}
\label{prop: Calabi inequality (Cone)}
$$
\int_X  \big|(d\,|u|^{\frac p2})\big|_g^2\ome^n\leq 
\frac{p^2n }{4(p-1)}\int_X u|u|^{p-2}(1-F)\ome^n,
$$
where $p>2$ with $p\del>2n-2$.
For simplicity, we denote by $K p$ the constant $\displaystyle{\frac{p^2n }{4(p-1)}}$.
\end{proposition}
\subsection{An inequality for the induction}
\bgn{lemma}\label{lem: A(p)}
Let $A(p)$ be a positive function of one variable $p \in \R$ which satisfies the following inequality, 
\bgn{equation}\label{eq: ineq A(p) 1}
A(p\e)^p \leq c_1 p\( A(p-1)\)^{p-1} 
\end{equation}
where $c_1$ is a constant and $\e=\frac n{n-1}>1$. 
We assume that there is a natural number $N>1$ such that $A(p)< c_0$ for all $p \in [N-1, N\e]$.
Then $A(p)$ satisfies the following inequality 
\bgn{equation}\label{eq:  ineq A(p)}
A(p)  \leq s_1 (s_2 p)^{ -\frac{s_3}p}, 
\end{equation}
where constants $s_1, s_2$ depend only on $c_0, c_1$ and $N$ and 
$s_3=\frac{2\e}{\e-1}=2n$.
In particular,  we have 
$$
\lim_{p\to\infty} A(p) =s_1
$$
which implies that $A(p)$ is bounded by a constant which does not depend on $p$.
\end{lemma}
\bgn{proof}
We take a constant $s_2$ which satisfies 
\bgn{equation}\label{eq:  s2}
2^{s_3}c_1p^{-1}{\e^{2(n-1)} }\leq {{s_2}^{2}},  
\end{equation}
for all $p>N$. ( It is possible because the left hand side is bounded.)
and then we choose a constant $s_1>1$ with 
\bgn{equation}\label{eq: s3}
c_0\leq s_1(s_2 p)^{-\frac{s_3}p},
\end{equation}
for all $p$
(It is also possible since $\lim_{p\to\infty}(s_2 p)^{-\frac{s_3}p}=1$.)
Then we shall show that $A(p)$ satisfies the inequality (\ref{eq:  ineq A(p)}) by the induction on $m$. 
The inequality (\ref{eq:  ineq A(p)}) holds for $p\in [N-1,N\e]$.
We assume that $A(p)$ satisfies the inequality (\ref{eq:  ineq A(p)}) for all $p \in [N-1, N\e^m]$. 
Then we shall show that the inequality holds for $p\e \in [N-1, N\e^{m+1}]$.
From our assumption, $A(p)$ and $A(p-1)$ satisfies the inequality (\ref{eq: ineq A(p)}) and 
then
applying (\ref{eq: ineq A(p) 1}), we have 
\bgn{align}
A(p\e)^p &\leq c_1 ps_1^{p-1}( s_2 (p-1))^{-s_3}
\leq (2^{s_3} c_1p)s_1^p(s_2p)^{-s_3}
\end{align}
by using $(\frac p{p-1})^{s_3}< 2^{s_3}$.
If we have the following inequality
\bgn{equation}\label{eq: 2.12}
(2^{s_3} c_1p)s_1^p(s_2p)^{-s_3}\leq s_1^p (s_2 p\e)^{-\frac{s_3}\e},
\end{equation}
then $A(p\e)$ satisfies the inequality (\ref{eq:  ineq A(p)}). 
Thus it suffices to show (\ref{eq: 2.12}).
The inequality (\ref{eq: 2.12}) is equivalent to 
\bgn{equation}\label{eq: 2.13}
(2^{s_3}c_1p) p^{-s_3+\frac{s_3}\e}\leq (\e)^{\frac{-s_3}\e} s_2^{\frac{s_3}n}
\end{equation}
Since $-s_3+\frac{s_3}\e=-\frac{s_3}n=-2$ and $\frac{s_3}\e=2(n-1)$,
(\ref{eq: 2.13}) is 
\bgn{equation}
2^{s_3}c_1p^{-1}\leq (\e)^{-2(n-1)} {s_2}^{2}
\end{equation}
This is the inequality which $s_2$ satisfies. 
Thus the result follows from the induction.
\end{proof}
\subsection{The openness}
Our proof of the theorem relies on the continuity method.
We consider a family of the \M-A equation parametrised by $s\in [0,1]$. 
For $F_s= 1-s+ sF$, we define 
\bgn{equation}\label{eq: s-Monge-Ampere}
\ome_{u_s}^n =F_s\ome^n
\end{equation}
Let $S$ be a subset of $[0,1]$ defined by 
$$
S=\{\, s\in [0,1]\, \,|\, \text{ \rm the equation (\ref{eq: s-Monge-Ampere}) has a solution } u_s\in C^{k+2,\a}_{\del}\,\},
$$
Note that $F-1\in C^{k,\a}_{\til\del}$, for all $\til\del\leq \del$.
Since $F_0=1$, the equation (\ref{eq: s-Monge-Ampere}) has the trivial solution $0$ for $s=0$. 
If we show that $S$ is open and closed, then there is a solution $u_1$ for $s=1$ which gives 
the solution of the \M-A equation in the theorem.
In order to show that $S$ is open, taking the derivative of the equation (\ref{eq: s-Monge-Ampere}) at $u_s$, 
we have the equation of the Laplacian 
\bgn{equation}\label{eq: eq of Laplacian s}
-\Psi_s \trian_s\dot{u_s}= \dot{F_s},
\end{equation}
where $\dot{u_s}\in C^{k+2,\a}_\del$ and $\dot{F_s}\in C^{k,\a}_{\del+2}$ and 
$\Psi_s$ is a bounded function given by 
$$
\Psi_s= \frac 12 \frac{\ome_{u_s}^n}{\ome^n}
$$
and $\trian_s$ is the Laplacian with respect to the K\"ahler metric $\ome_{u_s}$. 
Since $u_s\in C^{k+2,\a}_\del$, the metric $\ome_{u_s}$ is a conical K\"ahler metric. 
Then since $F-1\in C^{k,\a}_{2+\del}$,
it follows from the proposition \ref{prop: Laplacian(conical)} there is a unique solution $\dot{u_s}$ of the equation (\ref{eq: eq of Laplacian s}).
We define two Banach manifolds ${\cal M}=\{\, u\in C^{k+2,\a}_\del\, |\, \ome_u >0\, \}$.
${\cal N} =\{ \, f+1\, |\, f \in C^{k,\a}_{2+\del}\, \}$. 
Then the smooth map $\Psi  :{\cal M} \to {\cal N}$ is defined by 
$$
\Psi(u)=\frac{\ome_u^n}{\ome^n}.
$$
The differential of $\Psi$ at $u$ is the map $d\Psi_u : C^{k+2,\a}_\del \to C^{k,\a}_\del$ which is given by 
$$
d\Psi_u(\dot{u}) = -\Psi_u\trian_u,
$$ 
where $\trian_u$ is the Laplacian with respect to the K\"ahler metric $\ome_u$. 
Then from the proposition \ref{prop: Laplacian(conical)}, $d\Phi_u$ is isomorphism and it follows from the implicit function theorem that $\Phi$ is locally an isomorphism. 
Then $S$ is an open set.

\subsection{$C^0$-estimates}
We shall show the $C^0$-estimate of the \M-A equation 
$\ome_{u_s}^n = F_s\ome^n$. 
For simplicity, we write $u$ for $u_s$ and $F$ for $F_s$. 
We use $\ome^n$ as the (conical) volume form of the integrant in this subsection.
Recall that $u\in C^{k+2, \a}_\del$ and $F-1\in 
C^{k,\a}_{2+\del}$, for 
$0<\del <2n-2$. 
Substituting $f=|u|^\frac p2$ into the Sobolev inequality (proposition \ref{prop: Sobolev (conical)}), 
we have 
$$
\(\int_X |u|^{p\e}\ome^n \)^{\frac 1\e}\leq S^2 \int_X \big|\, d|u|^{\frac p2}\,|_g^2\,\ome^n
$$
Applying the inequality (proposition \ref{prop: Calabi inequality (Cone)}), we have 
\bgn{equation}\label{eq: 2.18}
\(\int_X |u|^{p\e}\ome^n \)^{\frac 1\e}\leq S^2Kp \int_X u|u|^{p-2}(1-F)\ome^n
\end{equation}
We set $A(p)$ by 
$$
A(p) = \(\, \int_X |u|^p\,\ome^n\,\)^{\frac 1p}
$$
and $c_1=S^2K\|1-F\|_{C^0}$ by the $C^0$-norm of $1-F$.
Then the inequality (\ref{eq: 2.18}) implies  
\bgn{equation}\label{eq: 2.19}
A(p\e)^p \leq c_1 p\( A(p-1)\)^{p-1}.
\end{equation}
Then we shall show an estimate of the $C^0$-norm by using the lemma 
\ref{lem: A(p)}. 
Recall $p\del >2n-2$ in the proposition \ref{prop: Calabi inequality (Cone)}. 
We set $r>1$ by $(p-1)r= p\e$ and define $s>1$ by 
$\frac 1s +\frac 1r=1$. Then it follows 
\bgn{equation}\label{eq: 2.20}
(\del +2)s >2n.
\end{equation}
Applying the H\"older inequality to the right hand side of (\ref{eq: 2.18}) with the exponents $r,s$, 
we have 
\bgn{equation}
A(p\e)^p\leq S^2Kp\( A( p\e)\)^{\frac{p\e}r} \, \( \int_X |1-F|^s\,\ome^n\,\)^\frac 1s
\end{equation}
Since $p-\frac{p\e}r=1$, we have 
$$
A(p\e)\leq S^2K p \( \int_X |1-F|^s\,\ome^n\,\)^\frac 1s.
$$
Note that from (\ref{eq: 2.20}), $\( \int_X |1-F|^s\,\ome^n\,\)<\infty$.
Thus for sufficiently large $N$, we see that 
$A(p)<c_0$ for all $p\in [N-1, N\e]$.
Then together with (\ref{eq: 2.19}), we can apply the lemma \ref{lem: A(p)} to obtain 
$$
A(p) =\(\, \int_X |u|^p\,\ome^n\,\)^{\frac 1p}
<C,
$$
where $C$ does not depend on $p$.
Thus we obtain an estimate of $C^0$-norm of $u$. 
\subsection{$C^2$-estimate}
There are many good references for an estimate of $C^2$-norm of 
a solution $u$ of the \M-A equation 
 \cite{Au}, \cite{Jo}, \cite{Na}, \cite{Ti}.
Bando-Kobayashi \cite{BK2} developed a way for the $C^2$-estimate
to apply the Schwartz lemma which
makes a comparison of two metrics $\ome$ and $\ome_u$ 
by using their curvature,
which is clear from a geometric point of view. 

We shall give a quick view of the $C^2$-estimate on our conical K\"ahler manifolds.
Let  
$\tr_{\ome_u}\ome$ be the trace of $\ome$ by the K\"ahler metric $\ome_u$. 
It follows that  
$\tr_{\ome_u}\ome>0  $ is the norm of $\ome$ by $\ome_u$ which gives the $C^2$-norm of $u$.
According to the Bando-Kobayashi's result, there is a constant $C$ depending only on 
the curvature of $\ome$ and the Ricci curvature of $\ome_u$ such that the following inequality holds,
$$-\square_u \log\tr_{\ome_u}\ome\geq -C(1+\tr_{\ome_u}\ome)
$$
where 
$\square_u$ is the complex Laplacian with respect to the K\"ahler metric 
$\ome_u$. 
Since $u$ is a solution of the \M-A equation, the Ricci curvature of $\ome_u$ is already bounded.
Then by using,
$$
-\square_u u =n-\tr_{\ome_u}\ome
$$
we have 
$$
-\square_u(\log\tr_{\ome_u}\ome-(C+1)u)\geq 
\tr_{\ome_u}\ome-(C+1)n-C
$$
Since the function $u$ decays exponentially with order $O(e^{-\del t})$,
the function $(\log\tr_{\ome_u}\ome-(C+1)u)$ is maximum at a point $x_0\in X$. Then
\bgn{align}
0\geq&-\square_u(\log\tr_{\ome_u}\ome-(C+1)u)(x_0)\\
\geq&
(\tr_{\ome_u}\ome -(C+1)n-C)(x_0)
\end{align}
Hence
$$
(C+1)n+C\geq \tr_{\ome_u}\ome(x_0)
$$
It follows that 
$$
(\log\tr_{\ome_u}\ome-(C+1)u)(x_0)\geq 
(\log\tr_{\ome_u}\ome-(C+1)u)(x), \quad \forall x\in X
$$
Thus we have the upper bound for $\tr_{\ome_u}\ome$,
$$
\sup_X\log\tr_{\ome_u}\ome <\log ((C+1)n+C)+2(C+1)\|u\|_{C^0}.
$$
By using the \M-A equation $\ome_u^n =F\ome^n$, 
we have $C^2$-estimate of $u$ and we obtain an apriori constant $C>0$ such that 
$$C^{-1}\ome< \ome_u< C\ome.$$

\subsection{$C^0_{\til\del}$-estimate}
In this subsection, we shall give an estimate of the weighted $C^0$-norm $u$ of a solution of the \M-A equation. 
In the process of our estimate, we need to choose a smaller weight $\til\del$ with $0<\til \del<\del$. 
We can obtain a way of the weighted $C^0$-estimate by applying
the method by Joyce,
(see section 8 in \cite{Jo}).
Our notation are the same as in the previous subsection.

At first we have the following weighted inequality,
\bgn{lemma} \label{weighted CY inequality}
 For a constant $q>0$ satisfying 
$$
q<p\del-2n+2
$$
we have 
\bgn{align*}
\int_X  r^q\big|(d\,|u|^{\frac p2})\big|_g^2\ome^n\leq &
\frac{p^2m }{4(p-1)}\int_X r^q u|u|^{p-2}(1-F)\ome^n\\
&+q^2C\int_Xr^{q-2} |u|^p\ome^n
\end{align*}
where $r=e^t$ and a constant $C$ depends only on 
the $C^2$-norm of $u$ given in the previous section.
\end{lemma}
\bgn{proof}
The condition 
$$
q<p\del-2n+2
$$
implies that the Stokes theorem holds,
$$
\int_X d\( r^qu|u|^{p-2}d^c u\w(\ome^{n-1}+\cdots+\ome_u^{n-1})\)=0.
$$ 
Then the result follows as in the way of the proposition \ref{prop: Calabi inequality (Cone)}
\end{proof}

The following lemma is a slight generalization of the lemma \ref{lem: A(p)}
\bgn{lemma}\label{lem: B(p)}
Let $B(p)$ be a positive function with one variable $p\in\R_{>0}$. 
We assume that  $B(p)$ satisfies the following,
 \bgn{align}
 B(p\e)^p\leq c_1p B(p-1)^{p-1}+(c_2p^2+c_3) B(p)^p
 \end{align}
 where  $\e=\frac{n-1}n$ and 
  there are a natural number $N$ and a positive $c_0$ such that 
  $B(p)<c_0$ for all
  $ p\in [N-1,N\e]$. 
  Then there are constants 
  $s_1, s_2$ depending only on $c_0,c_1,c_2,c_3$  such that for all $p>N$, $B(p)$ satisfies the inequality,
  $$
  B(p)\leq s_1( s_2p)^{-\frac{s_3}p}
  $$
 where $s_3=\frac{2\e}{\e-1}$.
 In particular, we have 
 $$
 \lim_{p\to\infty}B(p) =s_1
 $$
 Thus
 $B(p)$ is bounded by a constant which does not depend on $p$. 
\end{lemma}
\bgn{proof}
The proof is essentially same as the one in lemma \ref{lem: A(p)}.
\end{proof}

We shall start our 
$C^0_{\til\del}$-estimate, paying our attention on the weight.
We choose a weight $\til\del$ satisfying 
$0<\til\del<\del$.
and try to obtain a $C^{0}_{\til\del}$-norm of $u$.
We set 
$$q=p\til\del -2n+2$$
Then $q$ satisfies
\bgn{equation}\label{eq: stokes condition}
-p\del  +2n-2 +q <0
\end{equation}
which is the condition of the lemma \ref{weighted CY inequality}.
We define $B(p)$ by the weighted $L^p$-norm of $u$ with respect to the cylindrical metric.
$$
B(p) =\|u\|_{\dstyle L^p_{0,\til\del}}=\(\int_X |e^{\til\del t} u|^p \vol_{\cyl} \)^{\frac1p}
$$
Note $\ome^n =e^{2n t}\vol_{\cyl}$. Then we have 
\bgn{align}\label{eq: 25}
\int_X e^{\e qt}|u|^{p\e}\ome^n=&\int_X 
e^{\e p \til\del t }|u|^{p\e}  e^{(-2n +2)\e t} e^{2n t} \vol_{\cyl}
\end{align}\label{eq: weighted C0-estimate 25}
Since $\e=\frac n{n-1}$, we have  
$(-2n+2)\e+ 2n = 2(1-n)\frac n{n-1} + 2n =0$.
Then it follows
\bgn{align}\label{eq: weighted C0-estimate 26}
\int_X e^{\e qt}|u|^{p\e}\ome^n=&\int_X |e^{\til\del t}u|^{p\e}\vol_{\cyl}
=B(p\e)^{p\e}
\end{align}
We substitute 
$e^{{\frac q2}t } |u|^{\frac p2}$  into the Sobolev inequality in the
proposition \ref{prop: Sobolev (conical)}.
Since in the Sobolev inequality, the norm is given by the conical metric $\ome$ and the volume form is also 
  $\ome^n$, the left hand side of the inequality is 
 \bgn{align}
 \|e^{{\frac q 2}t}|u|^{\frac p2}\|^2_{\dstyle L^{2\e}}=&
 \(\int_X e^{q\e t}|u|^{p\e}\ome^n\)^{\frac 1\e}=B(p\e)^p
 \end{align}
 and using the norm by the cone metric, we find that the right hand side is 
 \bgn{align}
 \big\| d \(e^{\frac{q\til\del }2t}|u|^{\frac p2}\)\big\|^2_{\dstyle L^2}\leq&
 C\int_X e^{q t}\big|d|u|^{\frac q2}\big|^2\ome^n+
 Cq\til\del\int_X e^{(q-2) t}|u|^p\ome^n
 \end{align}
 The second term is given by 
 \bgn{equation}
 \int_X e^{(q-2)\til\del t}|u|^p\ome^n=\int_X |e^{\til\del t} u|^p
 e^{-2n t} e^{2n t} \vol_{\cyl}
 =B(p)^p
 \end{equation}
 Since  $e^{(2+\til\del )t}|1-F|<C$, we find that 
 \bgn{align}
 \int_X e^{q t}|u|^{p-1}|1-F|\ome^n =&\int_X e^{-\del t }e^{(q-2)t} |u|^{p-1}
 e^{(2+\til\del) t}|1-F|\ome^n \\
 \leq & C\int_X e^{-\til\del t}e^{p\til\del t}|u|^{p-1}\vol_{\cyl}\\
 \leq &C\int_X|e^{\til\del t }u|^{p-1}\vol_{\cyl}\\=&CB(p-1)^{p-1}
 \end{align}
Substituting these into the Sobolev inequality and 
 combining the inequality in the proposition \ref{weighted CY inequality}
 we obtain 
 \bgn{align}
\|u_s\|^p_{\dstyle L^{p\e}_{0,\til\del} }
<&C_1 p 
\|u_s\|^{p-1}_{\dstyle L^{p-1}_{0,\til\del}}+(C_2p^2+C_2)
\|u_s\|^p_{\dstyle L^p_{0,\til\del}}\\
\end{align}
The inequality implies 
\bgn{equation}
B(p\e)^p< C_1 p B(p-1)^{p-1}+(C_2p^2+C_2)B(p)^p
\end{equation}
We apply a simple trick to show a bound of $B(p)$ from $A(p)$ with respect to 
the weight.
Recall 
$$A(p) =\(\int_X |u|^p\ome^n \) ^\frac 1p=
\(\int_X |u|^p e^{2n t} \vol_{\cyl}\)^\frac 1p
$$
For $p\del >2n-2$, we already have a bound $C$ of 
$A(p)$ in the subsection of $C^0$-estimate.
On the other hand, 
$$
B(p) =\( \int_X  |u|^pe^{p\til\del t}\vol_{\cyl}\)^{\frac 1p}
$$
Thus 
for $p $ satisfying  $2n-2 < p\til\del <2n$,
we have  $B(p)< A(p)< C$, (where $\til\del <\del$.)
Applying the lemma \ref{lem: B(p)} to our $B(p)$, 
 We obtain the estimate  
$B(p)<C$, where $C$ does not depend on $p$. 
Thus it follows that the the $C^0$-norm of
$e^{\til\del t}u$ is bounded.

\subsection{$C^0_\del$-estimate by the maximum (minimum) principle}
We shall obtain an $C^0_\del$-estimate by using the maximum 
(minimum) principle.
\bgn{proposition}
There is a constant $C$ which does not depend on $s\in S$ such that $\|e^{\del t} u_s\|_{C^0}<C$.
\end{proposition}
\bgn{proof}
We shall use the equation (\ref {eq: eq of Laplacian s}) in subsection 2.5 to 
obtain an $C^0_\del$-estimate of a solution of \M-A equation,
\bgn{equation}
-\Psi_s \trian_s\dot{u_s}= \dot{F_s},
\end{equation}
where $\dot{u_s}=\frac{d}{dt}u_s$ and 
$\dot{F_s}=F-1$. 
Note that the family of solutions $\{u_s\}_{s\in S}$ is differentiable 
with respect to the parameter $s$ from the implicit function theorem. 
It suffices to show a  uniform $C^0_\del$ estimate of 
$\dot{u_s}$ in order to obtain a uniform $C^0$-estimate of $u_s$
for $s\in S$.
We already have the $C^0_{\til\del}$-estimate.
Note that $\Psi_s$ is bounded.  Then as in compact K\"ahler manifolds, we have $C^{k,\a}_{\til\del}$-estimate, for $0<\til\del<\del$.
It implies that we have the $C^{k,\a}_{\til\del}$-estimate of a solution 
$\dot{u}_s$ of the equation (\ref{eq: eq of Laplacian s}) which does not depend on $s$.
The difference of every coefficients of the Laplacians $\trian_s -\trian$ 
decay exponentially with the order $e^{-(2+\til\del)t}$.
Then we have constants $C_3$ and $T>0$ such that 
$$
\trian_se^{-\del t}\geq C_3 e^{-(2+\del)t},
$$
on the region $\{ t> T\}$.
We take a constant $C$ which satisfies the two inequalities, (\ref{eq: Maximum principle-2}) and (\ref{eq: Maximum principle-3}), 
where $\|\dot{u}_s\|_{C^0}< C_1$ and $\|e^{(2+\del)t}\Psi_s^{-1}\dot{F}_s\|_{C^0}<C_2$.
Then we can apply the same method as in lemma \ref{lem: maximum principle}
to obtain 
$$
\|e^{\del t}\dot{u}_s\|_{C^0}<C,
$$
where $C$ does not depend on $s$. 
Hence we have a $C^0$-estimate of $e^{\del t}\dot{u_s}$. Thus
we obtain a $C^0_{\del}$-estimate of $u_s$.
\end{proof}
\bgn{proof}{\it of the existence theorem \ref{th: existence CY on conical}}
\,We obtain a $C^{2,\a}$-estimate of $u$ by applying the general method of 
the $2$-nd order elliptic differential equations to our conical K\"ahler manifolds  as in the case of compact K\"ahler manifolds.
This method was developed in \cite{Ev}, \cite {Kry}. 
(for instance, see \cite{Tr}).
Successively we apply the Schauder estimate to obtain the $C^{k+2,\a}$-estimate of $u$ 
(see Theorem 6.2 and Theorem 17.15 in\cite{G-T}).
This procedure is explained in page 89 \cite{Siu}.
 We have the equation, 
 \bgn{equation}\label{eq: weighted Schauder}
 \trian_se^{\del t}\dot{u}_s =H
 \end{equation}
 where the $C^{k,\a}$-norm of $H$ is bounded.
 The $C^{k,\a}$-norms of coefficients of the Laplacian $\trian_s$ are bounded.
 We already have the bound of $C^0_\del$-norm of $u$.
 Then applying the Schauder estimate to the equation (\ref{eq: weighted Schauder}), 
 we obtain an estimate of the weighted norms $C^{k+2,\a}_\del$ of $u$. 
Hence it follows that the set $S$ is closed by the Ascoli-Arzel\`a lemma. 
Since $S$ is open, it follows that $S=[0,1]$. 
Thus we have a solution of the \M-A which gives the 
Ricci-flat conical K\"ahler metric on $X$.
\end{proof}\par
\medskip
\noindent
\bgn{proof}{\it  of  the theorem \ref{th: uniqueness theorem } }
Since $\ome'=\ome_u= \ome+\sqrt{-1}\pa\ol\pa u$, 
and  both $\ome$ and $\ome'$ satisfy the \M-A equation,
then we have 
$\Ome\w\ol\Ome = c_n \ome^n =c_n \ome_u^n$. 
It follows that $F\equiv 1$. 
Then  apply the inequality of the proposition \ref {prop: Calabi inequality (Cone)}
for $\del p>2n-2$, we have 
$$
\int_X|(d|u|^{\frac p2}) |^2_g \ome^n =0.
$$
Hence $u\equiv$ constant and if follows $\ome=\ome'$.
\end{proof}\par\medskip

\noindent
\bgn{proof}{\it of theorem }\ref{th: uniqueness in class}\,
We set $\a=\ome-\ome'$. Then the $d$-exact $2$-from $\a$ decays with the order $O(e^{-\lam t})$, i.e., $\a\in C^k_\lam$ with respect to the cylindrical metric.
Let $f$ be the real function given by the contraction $\w_\ome\a$ in terms of the K\"ahler form $\ome$, 
where the real function $f$ lies in $C^k_{\lam+2}$. 
For $n-2<\lam< 2n-2$, we put $\lam=\lam_1$. We have a solution $u\in C^{k+2}_{\lam_1}$ such that $\frac12\trian_\ome u =f$, where $\trian_\ome$ denotes the Laplacian with respect to $\ome$. 
In the cases of $\lam\geq 2n-2$, then we also have a function $u\in C^{k+2}_{\lam_1}$ such that 
$\frac12\trian_\ome u =f$ for $n-2<\lam_1<2n-2$ since $f\in C^k_\lam\subset C^k_{\lam_1}$ for $\lam_1<\lam$.
Put $\b =\a-\frac1{\sqrt{-1}}\pa\ol\pa u$. Then using $\frac1{\sqrt{-1}}\w_\ome\pa\ol\pa u =\square_\ome u =\frac12\trian_\ome u$, 
we have $\w_\ome \b=0$. Thus $\b$ is a real primitive form of type $(1,1)$ and we have 
\bgn{equation}\label{eq:proof of theorem 1.8-1}
\b\w\b\w\frac{\ome^{n-2}}{(n-2)!}=-|\b|^2\frac{\ome^n}{n!}
\end{equation}
We denote by $H^\bullet_\lam(X)$ the cohomology group which are defined by 
$$
H^p_\lam(X) =\frac{  \{\, a\in C^k_\lam(\w^p)\, |\, da=0\, \}  }{ \{db\, |\, b\in C^{k+1}_\lam(\w^{p-1})\, \}},
$$
where $C^k_\lam(\w^p)$ denotes $p$-forms which decay with the order $O(e^{-\lam t})$.
Then it turns out that $H^p_\lam (X)\cong H^p_{\cpt}(X)$ for $\lam >0$, where $H^p_{\cpt}(X)$ is the compactly supported cohomology group.
Since $H^1(S)=\{0\}$, we have the exact sequence, 
$$
0\to H^2_{\cpt}(X)\to H^2(X)\to H^2(S)\to \cdots.
$$
Hence the map $H^2_{\cpt}(X)\cong H^2_\lam(X)\to H^2(X)$ is injective for $\lam>0$.
Since $\b$ is $d$-exact which lies in $C^k_{\lam_1}(\w^2)$ for $n-2<\lam_1$, it follows that 
there is a $1$-from $\eta\in C^{k+1}_{\lam_1}(\w^1)$ such that $d\eta=\b$.
Thus we have 
$\b\w\b\w\frac{\ome^{n-2}}{(n-2)!}= d\eta\w\b\w\frac{\ome^{n-2}}{(n-2)!}.$
Since $\ome$ lies in $C^\infty_{-2}$, we see that $\eta\w\b\w\ome^{n-2}$ decays with the order $O(r^{-2(\lam_1-n+2)})$
for $\lam_1>n-2$. 
Then it follows from the Stokes theorem that we have 
$$
\int_X \b\w\b\w\frac{\ome^{n-2}}{(n-2)!}=\int_X d\eta\w\b\w\frac{\ome^{n-2}}{(n-2)!}=0
$$
Then from (\ref{eq:proof of theorem 1.8-1}), we see that $\b=0$. 
Thus $\a=\frac1{\sqrt{-1}}\pa\ol\pa u$ for $u\in C^{k+2}_{\lam_1}$. 
Since $\lam_1>0$, the result follows from the theorem \ref{th: uniqueness theorem }.
\end{proof}


\section{Sasakian structures and K\"ahler structures on the cone}
We will give a brief explanation of Sasakian manifolds from a view point of 
a correspondence to K\"ahler structures on the cones.
A notion of Sasaki manifolds was introduced by Sasaki and now
there is  a good reference 
on Sasaki geometry\cite{BG} in which the material in this section 
can be found. 
Let $S$ be a compact manifold of dimension $2n-1$.
Note that $S$ does not have a boundary.
The cone of $S$ is the product $\Bbb R_{>0}\times S$ with $r=e^t\in \R_{>0}$.
A Riemannian metric  $g_{\ss S}$ on $S$ 
yields the cone metric $g$ on $C(S)$ by
\bgn{equation}\label{eq: cone metric}
g={dr}^2 +r^2 g_{\ss S}
\end{equation}
In this section a metric on $C(S)$ is always the cone metric $g$ which is given by a metric $g_S$ on $S$
\bgn{definition}\label{def: Sasakian manifolds}
A $(2n-1)$ dimensional Riemannian manifold
$(S, g_{\ss S})$ is  a Sasakian manifold if there is a complex structure $J$ on the cone such that $(C(S), g, J)$ is a K\"ahler manifold. 
\end{definition}
This is a relevant definition of Sasakian manifolds focusing on the relation to 
K\"ahler geometry, however which is different from the ordinary definition.

We shall explain the correspondence to the ordinary one in which  
interesting geometric structures are included.
The K\"ahler structure $(C(S), g, J, \ome)$ as in definition \ref{def: Sasakian manifolds}
gives geometric structures on $S$.
At first we regard $S$ 
as the hypersurface  $\{x\in C(S)\, |\, r(x) =1\, \}$ in $C(S)$.
We denote by $r\frac{\pa}{\pa r}=\frac\pa{\pa t}$ the vector field on $C(S)$ which is defined by the translation of the cone in terms of the cylinder parameter
$t$.
The following lemma is known, for instance, see appendix 
\cite{MS}.
\bgn{lemma}
Let $(C(S), g, J, \ome)$ be the K\"ahler manifold of a Sasakian manifold $S$. 
Then we have $\L_{\frac{\pa}{\pa t}}J=0$
\end{lemma}
The Lie derivative  $\L_{r\frac{\pa}{\pa r}}J=0$ implies that  $J$ is invariant under the translations with respect to $t$, 
in other words, the vector field $\frac{\pa}{\pa t}$ is the real part of a holomorphic vector field.
\bgn{proof}
Let $\nab$ be the Levi-Civita connection of the cone metric $g$. 
Then we see that 
$\nab_u\frac{\pa}{\pa t}= u$, for all vector field $u$. 
For any vector fields $u, v$, we have 
\bgn{align*}
(\L_u J)v =&[u, Jv]-J[u,v]\\
=&\nab_u(Jv)-\nab_{Jv}u-J\nab_u v+J\nab_v u\\
=&(\nab_u J)v-\nab_{Jv}u +J\nab_v u
\end{align*}
Since $g$ is a K\"ahler metric, we have $\nab J=0$. 
Then substituting $u=\frac{\pa}{\pa t}$, we obtain 
$$
(\L_{\frac{\pa}{\pa t}}J)v =-\nab_{Jv}\frac{\pa}{\pa t}
+J\nab_v \frac{\pa}{\pa t} =0.
$$
\end{proof}
We define a vector field $\xi=Jr\frac{\pa}{\pa r}=J\frac\pa{\pa t}$ on 
$C(S)$. Since $g$ is a Hermitian metric, $\xi$ is orthogonal to $\frac{\pa}{\pa t}$ which implies that 
$\xi$ is a vector field along the hypersurface $S$. 
The restriction of  $\xi$ to $S$ is denoted by 
$$\xi_{\ss S}=\xi|_{\ss S}$$
We have the complex structure $J^*$ which acts on $1$-forms
by $(J^*\theta)(v) = \theta( Jv)$, for $\theta\in T^*C(S)$ and
$v\in TC(S)$. 
Then a $1$-form $\eta$ on $C(S)$ is given by 
$$\eta=- J^*\frac {dr}r=-J^* dt$$
and we denote by  $\eta_{\ss S}$ 
the restriction of $\eta$ to $S$.
(For simplicity, we write $J$ for $J^*$ from now on.)
Then it follows from the definition that 
$\eta(\xi) =\eta_{\ss S}(\xi_{\ss S}) =1$.
The K\"ahler form $\ome$ is defined by 
$\ome(u,v) =g(Ju, v)$, 
for $u, v\in TC(S)$.  
The Lie derivative  $\L_{\frac {\pa}{\pa t}}$ is induced from the translation of $t$-direction on $C(S)$. The group of one parameter transformations $f_\lam$ on $C(S)$ is given by
 $$f_\lam( r, y) = ( \lam r, y),\qquad y\in S
 $$ and 
 $$\frac d{d\lam}f_\lam^*|_{\lam=0}=\L_{r\frac \pa{\pa r}}
 $$
 
 From (\ref{eq: cone metric}), we have 
 $f_\lam^* g= \lam^2 g$ and $\L_{r\frac \pa{\pa r}}g= 2g$.
Since $\L_{r\frac \pa{\pa r}}J=0$, we have 
$\L_{r\frac \pa{\pa r}}\ome=2\ome$ .
Since $d\ome=0$, applying the formula of the Lie derivative 
$\L_u =i_u d+ d i_u$, we have 
\bgn{equation}\label{eq: 2ome}
d i_{r\frac \pa{\pa r}}\ome= 2\ome
\end{equation}
Then we have
\bgn{lemma}\label{lem: eta}
$$
i_{r\frac{\pa}{\pa r}}\ome=i_\xi g =r^2\eta
$$
\end{lemma}
\bgn{proof}
Since $\xi =J\frac{\pa}{\pa t}$, 
and $\eta=-J^*dt$, we have $\eta(\xi)=1$.
Then $i_\xi g (\xi ) =r^2\eta(\xi)=r^2$.
Since $J\frac{\pa}{\pa t} \in TS$, we obtain $\eta(\frac{\pa}{\pa t}) = -dt(J\frac{\pa}{\pa t}) =0$.
Let $\lan \frac{\pa}{\pa t} , J\frac{\pa}{\pa t}\ran^\perp$ 
be the orthogonal subspace to the space spanned by two vector fields $\frac{\pa}{\pa t} $ and $ J\frac{\pa}{\pa t}$ . Then $\lan \frac{\pa}{\pa t} , J\frac{\pa}{\pa t}\ran^\perp$ is invariant under the action of $J$ and  we have   $\eta( u ) =-dt( Ju)=0$ for  $ u \in 
\lan \frac{\pa}{\pa t} , J\frac{\pa}{\pa t}\ran^\perp$.
Thus $i_\xi g =r^2\eta$
\end{proof}\par
Applying 
lemma \ref{lem: eta} to the equation
 (\ref{eq: 2ome}), 
since $\eta=-J^*\frac{dr}r$, we have 
\bgn{align}
2\ome=& dr^2\eta =-d( J^*r dr)=-\frac12 dJ^* dr^2
\end{align}
Hence 
$$
\ome= \frac{\sqrt{-1}}2\pa \ol\pa r^2
$$

Since $2\ome =2rdr\w\eta+r^2d\eta$ is a symplectic structure on $C(S)$,
$$
(2\ome)^n=2nr^{2n}dt\w\eta\w(d\eta)^{n-1}\neq 0
$$
It implies that
$\eta_{\ss S}\w(d\eta_{\ss S})^{n-1}\neq 0$. We see that
$\eta_{\ss S}$ is a contact structure on  $S$ and
$\xi_{\ss S}$ is the Reeb vector field.
We define a distribution $D$ of dimension $2n-2$ by 
$D=\ker \eta_{\ss S}= \{\, u\in TS\, |\, \eta_{\ss S}(u) =0\, \}$,  
then $D= \lan \frac{\pa}{\pa t} , J\frac{\pa}{\pa t}\ran^\perp$
and $D$ is invariant under the action of $J$.
We define a section 
$\Phi_{\ss S}\in$End$(TS)$ by
$$
\Phi_{\ss S}(v) =
\bgn{cases}
&Jv \qquad  (v \in D)\\
&0\qquad \,\,\,\,(v = \xi )
\end{cases}
$$
Together with the contact structure 
 $\eta_{\ss S}$ and the Reeb vector field $\xi_{\ss S}$ on $S$,
$\Phi_{\ss S}\in$End $(TS)$ gives the Riemannian metric on $S$ by 
\bgn{equation}\label{eq: g S}
g_{\ss S}(u, v) =\eta_{\ss S}\otimes \eta_{\ss S}( u, v) +d\eta_{\ss S}( u, \Phi_{\ss S} v),
\end{equation}
for $u,v \in TS$. 
By the same procedure, an almost K\"ahler structure on $C(S)$ as in definition \ref{def: Sasakian manifolds} gives the structure $(\eta_{\ss S}, \xi_{\ss S}, \Phi_{\ss S})$ on $S$.
When $J$ is integrable, the corresponding structure on $S$ admits suitable properties. 
\bgn{lemma}
If a almost complex structure $J$ on $C(S)$ is integrable,
we have 
$$
\L_\xi J=J\L_{\frac{\pa}{\pa t}}J 
$$
\end{lemma}
\bgn{proof}
Since $J$ is integrable, the Nijenhuis tensor vanishes,
\bgn{equation}
\label{eq: Nijenhuis}
[J\frac{\pa}{\pa t}, Ju ]=J[J\frac{\pa}{\pa t}, u]+J[\frac{\pa}{\pa t},Ju]
+[\frac{\pa}{\pa t}, u]
\end{equation}
where $u \in TC(S)$.  Since $\xi=J\frac{\pa}{\pa t}$, we have  
\bgn{align}
\L_\xi J(u) =& [J\frac{\pa}{\pa t}, Ju]-J[J\frac{\pa}{\pa t}, u]\\
\L_{\frac{\pa}{\pa t}}J(u) =&[\frac{\pa}{\pa t}, Ju]-J[\frac{\pa}{\pa t},u]
\end{align}
Thus from 
(\ref{eq: Nijenhuis}), we have $\L_\xi J=J\L_{\frac{\pa}{\pa t}}J $.
\end{proof}\par
Then $\frac{\pa}{\pa t}-\sqrt{-1}J\frac{\pa}{\pa t}=
\frac{\pa}{\pa t}-\sqrt{-1}\xi$ is a holomorphic vector field on $C(S)$.
In particular, $[\frac{\pa}{\pa t}, \xi]=[\frac{\pa}{\pa t}, J\frac{\pa}{\pa t}]=0$.
\bgn{lemma}
If $\L_{\xi}J=0$, then  $\L_{\xi_{\ss S}}\eta_{\ss S}=0$.
\end{lemma}
\bgn{proof}
Since $\L_\xi \eta = -\L_\xi Jdt=-J\L_{\frac{\pa}{\pa t}}dt= 0$, 
the result follows from   
$i_{\ss S}^*\L_\xi \eta=\L_{\xi_{\ss S}}\eta_{\ss S}$. 
\end{proof}

Thus if $J$ is integrable, then $\L_\xi \eta=0$.
This also implies that 
$i_\xi d\eta =0$ and $d\eta_{\ss S}$ is a basic form on $S$.
Since $\L_{\xi_{\ss S}}\Phi_{\ss S}=0$ , (\ref{eq: g S}) yields $\L_{\xi_{\ss S}} g_{\ss S}=0$ and we see that 
$\xi_{\ss S}$ is a Killing vector filed on $(S, g_{\ss S})$. 

Hence we obtain the structure $(\eta_{\ss S}, \xi_{\ss S}, \Phi_{\ss S})$ on 
$S$ from a K\"ahler structure with $\L_{\frac{\pa}{\pa t} }J=0$ and a cone metric. 
Conversely we shall construct the K\"ahler structure as in definition \ref {def: Sasakian manifolds} on $C(S)$ from the structure  $(\eta_{\ss S}, \xi_{\ss S}, \Phi_{\ss S})$ on $S$.
 Let $\eta_{\ss S}$ be a contact structure on $S$.
Since the Darboux's theorem gives the standard form of a contact structure,
it follows that  
$D=\ker\eta_{\ss S}$ is a $2n-2$ dimensional distribution on which 
$d\eta_{\ss S}$ is a non-degenerate $2$-form.
 
Then a section $\Phi_{\ss S}\in$End $(TS)$ is defined as an almost complex structure on $D$ with 
$\Phi_{\ss S}(\xi_{\ss S})=0$.
Such a section $\Phi_{\ss S}$ is compatible with $\eta_{\ss S}$ if the following condition hold
\bgn{itemize}
\item $d\eta_{\ss S}(\Phi_{\ss S} u , \Phi_{\ss S} v) =d\eta_{\ss S}(u, v) , \quad u, v \in D$
\item $d\eta_{\ss S} (u, \Phi_{\ss S} u) >0, \qquad ( u\neq 0 \in D).$
\end{itemize} 
These conditions imply that  a pair $(d\eta_{\ss S}, \Phi_{\ss S})$ is a Hermitian structure on $D$. 
Thus a compatible  pair $(\eta_{\ss S}, \Phi_{\ss S})$ gives 
a Riemannian metric $g_{\ss S}$ by (\ref{eq: g S}) on $S$.
Then we have the cone metric $g$ on the cone 
 $C(S)=\R_{>0}\times S$ by 
$$
g= dr^2+ r^2 g_{\ss S}
$$ 
 The tangent bundle  $TC(S)$ is decomposed into $ T\R\times TS$ 
 and we regarded $\xi_{\ss S}\in TS$ as the vector field $\xi$ on
$C(S)$. 
An almost complex structure $J$ on $C(S)$ is given by 
$$
J(r\frac{\pa}{\pa r} ) =\xi , \qquad J|_D = \Phi_{\ss S}
$$
Then $g$ is a Hermitian metric with respect to $J$ and the corresponding $2$-form
 $\ome$ is a symplectic form,
$$
2\ome= d(r^2\eta_{\ss S}) = 2rdr\w\eta_{\ss S}+r^2d\eta_{\ss S}
$$
that is, $(g, J,  \ome)$ is an almost K\"ahler structure on $C(S)$.
Since $\Phi_{\ss S}$ is  a section of  End($TS$), 
the induced $J$ is invariant under the translation with respect to $t$.
Thus 
$\L_{\frac{\pa}{\pa t}}J=0$.
Hence we obtain the almost K\"ahler structure as in 
definition \ref{def: Sasakian manifolds}
from the structure $(\eta_{\ss S} ,\xi_{\ss S}, \Phi_{\ss S})$.
\bgn{definition}
Let $S$ be a manifold of dimension $2n-1$ which admits 
 a contact structure $\eta_{\ss S}$ and a compatible structure $\Phi_{\ss S}$ to $\eta_{\ss S}$. 
A compatible pair $(\eta_{\ss S}, \Phi_{\ss S})$ is a Sasakian structure on $S$ if the corresponding almost K\"ahler structure $(C(S), g, J,\ome)$ is K\"ahlerian, that is, $J$ is integrable.
 \end{definition}
Hence our argument is reduced to the following, 
\bgn{proposition}\label{prop: Sasakian and kahler cone metric}
There is a one to one correspondence between Sasakian structures $(\eta_{\ss S}, \xi_{\ss S}, \Phi_{\ss S})$ on $S$ and 
K\"ahler structures $(g, J, \ome)$ on $C(S)$ consisting of a cone metric $g$ and a
translation-invariant complex structure $J$.
\end{proposition}

\section{Einstein-Sasakian structures and weighted Calabi-Yau structures}
As in section $1$, let $S$ be a compact manifold of dimension $2n-1$  and 
$C(S)$ the cone $\R_{>0}\times S$ with $r=e^t\in \R_{>0}$. 
We assume that $S$ is simply connected in this section.
Let 
$(C(S), J, g, \ome)$ be a K\"ahler structure on $C(S)$ corresponding to a Sasakian structure on $S$ by the proposition \ref {prop: Sasakian and kahler cone metric}.
Then as in section $1$, we have 
\bgn{align}
\ome= &\frac {\sqrt{-1}}2\pa\ol\pa r^2\\
g=&dr^2+ r^2 g_{\ss S}\\
\L_{\frac{\pa}{\pa t}}J&=0
\end{align}
We also see in section $1$
\bgn{equation}\label{eq: ome symplectic}
2\ome= d(r^2\eta_{\ss S}) = 2rdr\w\eta_{\ss S}+r^2d\eta_{\ss S}
\end{equation}
We assume that the canonical line bundle of $(C(S), J)$ is trivial and there is a nowhere-vanishing holomorphic $n$-form
$\Ome $. \\
We call an $n$-form $\Ome $ is of weight $n$ if $\Ome$ satisfies
\bgn{equation}
\L_{\frac{\pa}{\pa t}}\Ome =n\Ome
\end{equation}
Note that $\ome$ and $g$ satisfy 
\bgn{align}
&\L_{\frac{\pa}{\pa t}}\ome= 2\ome\\
&\L_{\frac{\pa}{\pa t}}g=2g
\end{align}
and in this sense these are of weight $2$.
Note that a pair $(\Ome, \ome)$ consisting of a holomorphic $n$-form and a K\"ahler form induces a K\"ahler structure $(J,\ome)$ on $C(S)$.
\bgn{definition}\label{def: weighted Calabi-Yau}
A pair consisting of a $d$-closed, holomorphic $n$-form $\Ome$ of weight $n$ and a K\"ahler form $\ome$ on $C(S)$ is a weighted Calabi-Yau structure 
if the induced K\"ahler structure $(J,\ome)$ corresponds to a Sasakian structure on $S$ as in proposition \ref{prop: Sasakian and kahler cone metric}
and satisfies the \M-A equation,
\bgn{equation}
 \Ome\w\sig(\ol\Ome)= \frac{(2\sqrt{-1})^n}{n!}\ome^n,
 \end{equation}
 where $\ol\Ome$ is the complex conjugate of $\Ome$ and 
 $\sig(\ol\Ome)$ is given by 
 $$
 \sig({\ol\Ome})=
 \bgn{cases}
 &+\ol\Ome, \qquad  n\equiv 0,1\, (\text{\rm mod }4),\\
  &-\ol\Ome,\qquad  n\equiv 2,3\, (\text{\rm mod }4).
 \end{cases}
 $$
 \end{definition}
We regard $S$ as the hypersurface 
$\{ \, t=0\,\}=\{\, r=1\,\}$ as before.
We denote by $r\frac{\pa}{\pa r}=\frac\pa{\pa t}$ the vector field on $C(S)$ which is defined by the translation of the cone in terms of the cylinder parameter
$t$.
The next lemma is known for instance, see appendix 

The Lie derivative  $\L_{r\frac{\pa}{\pa r}}J=0$ implies that  $J$ is invariant under the translations with respect to $t$, 
in other words, the vector field $\frac{\pa}{\pa t}$ is the real part of a holomorphic vector field.

 Let $i_{\ss S} : S\to C(S)$ be the inclusion of $S$ into $C(S)$ and
 $p_{\ss S} : C(S) \to S$ the projection to $S$.
 For a weighted Calabi-Yau structure  $(\Ome, \ome)$, using the interior product by the vector field $\frac{\pa}{\pa t}$, we define 
\bgn{align}\label{ali: psi, eta}
 &\psi=i_{\frac{\pa}{\pa t}}\Ome ,  \\
 &r^2\eta=i_{\frac{\pa}{\pa t}}\ome,    
 \end{align}
 and  restricting to $S$, we have  $\psi_{\ss S}:=i_{\ss S}^*\psi$、 $\eta_{\ss S}=i_{\ss S}^*\eta$.
Since the vector field $\frac{\pa}{\pa t}$ is the real part of a holomorphic vector field, it follows that $\psi$ is a holomorphic $(n-1)$-form on $C(S)$, 
however $\psi$ is not $d$-closed.
Since $\Ome$ is a $d$-closed form of weight $n$, we have
\bgn{equation}
\L_{\frac{\pa}{\pa t}}\Ome=di_{\frac{\pa}{\pa t}}\Ome =d\psi =n \Ome
\end{equation}
Since $(dt-\sqrt{-1}\eta)=dt+\sqrt{-1}Jdt$ is a form of type $(1,0)$
it follows from (\ref{ali: psi, eta}) that 
\bgn{equation}\label{eq: Ome=}
\Ome =(dt-\sqrt{-1}\eta)\w\psi
\end{equation} 
Restricting to $S=\{ t=0\}$, we have 
$$
d\psi_{\ss S}= -n\sqrt{-1}\eta_{\ss S}\w \psi_{\ss S}
$$
Since $(\Ome, \ome)$ 
satisfies the \M-A equation, substituting 
(\ref{eq: Ome=}), we have 
\bgn{align*}
\Ome\w\sig(\ol\Ome)=&(dt-\sqrt{-1}\eta)\w\psi \w\sig(\ol\psi)\w (dt+\sqrt{-1}\eta)\\
=&2\sqrt{-1}dt\w\eta\w\psi\w\sig(\ol\psi)
\end{align*}
\text{\rm it also follows from (\ref{eq: ome symplectic}) } \\
\bgn{align*}
\frac{(2\sqrt{-1})^n}{n!}\ome^n=&\frac{(\sqrt{-1})^n}{n!} 2n r^{2n}dt\w
\eta \w(d\eta)^{n-1}
\end{align*}
Then by applying the interior product $i_{\frac{\pa}{\pa t}}$ to 
both sides of the \M-A equation,
we have 
\bgn{align}\label{ali: r-Sasaki-M-A}
\eta\w\psi\w\sig(\ol{\psi})=
c_{n-1}r^{2n}
\eta\w(\frac 12d\eta)^{n-1}
\end{align}
By restricting to 
$S=\{ \, r=1\,\}$, we have 
\bgn{align}\label{ali: r-Sasaki-M-A on S}
\eta_{\ss S}\w\psi_{\ss S}\w\sig(\ol{\psi}_{\ss S})=
c_{n-1}
\eta_{\ss S}\w(\frac 12d\eta_{\ss S})^{n-1}
\end{align}
Let $D$ be the subbundle of $TS$ given by the kernel of $d\eta_{\ss S}$.
Then $\psi_{\ss S}$ defines 
a section $\Phi_{\ss S}$ on End$\,(D)$ which is an almost complex structure and 
then $\psi_{\ss S}$ is a nowhere vanishing section of the canonical line bundle 
$K_D$ of $D$.
Since $\ome$ is K\"ahlerian, we see that $d\eta_{\ss S}$ is a Hermitian form 
on $D$ and then $(\eta_{\ss S}, \Phi_{\ss S})$ is an Sasakian structure on $S$.
An Einstein-Sasakian structure on $S$ is a Sasakian structure whose metric $g$ is an Einstein metric on $S$.
\bgn{proposition}\label{prop: Einstein-Sasakian and weighted Calabi-Yau}
Let $(\Ome,\ome)$ be a weighted Calabi-Yau structure on $C(S)$. 
Then $(\Ome,\ome)$ corresponds to an Einstein-Sasakian structure on $S$ under the correspondence in section $1$.
Conversely, an Einstein-Sasakian structure on $S$ gives 
a weighted Calabi-Yau structure $(\Ome,\ome)$ on $C(S)$, 
that is, there is a one to one correspondence between 
Einstein-Sasakian structures on $S$ and weighted Calabi-Yau structures on the cone $C(S)$.
  \end{proposition}
 As though a correspondence as in the proposition is already known among experts in Sasakian geometry (see \cite{BG} for instance),
our point of view may be rather relevant since we focus on weighted Calabi-Yau structures on $C(S)$ which exactly correspond to Einstein-Sasakian structures on $S$.  
  In the appendix we shall give a proof. 
 
\section{Ricci-flat conical K\"ahler metrics on crepant resolutions of normal isolated singularities }
Let $X_0$ be an affine variety of complex dimension $n$ with only normal isolated singularity at $p\in X_0$. 
In this section, we assume that the complement $X_0\bsh\{p\}$ is biholomorphic to the cone 
$C(S)$ of an Einstein-Sasakian manifold $S$ of real dimension $2n-1$. 
As we see in section $3$ and $4$,  there is the weighted Calabi-Yau structure $(\Ome, \ome_0)$ on the cone $C(S)$ and the Ricci-flat K\"ahler cone metric  $\ome_0$ which is 
given by 
\bgn{equation}
\ome_0=\frac{\sqrt{-1}}2 \pa\ol\pa r^2
\end{equation}
\bgn{theorem}\label{th: existence theorem crepant}
Let $X_0$ be an affine variety with only normal isolated singularity at $p$.
We assume that the complement $X_0\bsh\{p\}$ is biholomorphic to the cone 
$C(S)$ of an Einstein-Sasakian manifold $S$ of real dimension $2n-1$.
If there is a resolution of singularity $\pi : X\to X_0$ with trivial canonical line bundle $K_X$, 
then there is a Ricci-flat conical K\"ahler metric for every K\"ahler class of $X$. 
\end{theorem}
\bgn{remark}
Van Coevering \text{\rm\cite{Coe}} showed that there is a Ricci-flat conical K\"ahler metric in
 the K\"ahler class
which belongs to the compactly supported cohomology group $H_{\cpt}^2(X)$
 of $X$. 
 Our theorem shows that this kind of restricted condition is not necessary and 
 implies that the conjecture on  the existence of complete Ricci-flat K\"ahler metrics in \text{\rm\cite{MS}} is affirmative.
\end{remark}
The Ricci-flat K\"ahler cone metric $\ome_0$ is written as 
\bgn{equation}
\ome_0 = e^{2t}( 2dt\w d^ct + dd^c t),
\end{equation}
where $dd^ct$ is the transversal K\"ahler metric on the Sasakian manifold $S$ and the Reeb vector field $\xi_{\ss S}$
is given by $J\frac{\pa}{\pa t}$ which gives the Reeb foliation on $S$. 
Let $H^p_B(S)$ be the basic cohomology group  on $S$ with respect to the Reeb foliation. 
There are  the Hodge decomposition of the basic cohomology groups  and the Lefschetz decomposition as in 
K\"ahler geometry. 
Thus the second basic cohomology group $H^2_B(S)$ is decomposed into the followings,
\bgn{align}
&H^2_B(S,\C) = H^{1,1}_B(S)\oplus H^{2,0}_B(S)\oplus H^{0,2}_B(S),\\
&H^{1,1}_B(S,\R) = \R dd^c t\oplus P^{1,1}_B(S,\R),
\end{align} 
where $H^{p,q}_B(S)$ denote the basic Dolbeault cohomology group of type $(p,q)$ and 
$\R dd^ct$ is the one dimensional space 
generated by the transversal K\"ahler form $dd^c t$
and $P^{1,1}_B(S)$ is the basic primitive cohomology group of type $(1,1)$ with respect to $dd^ct$.
Since the transversal complex structure on $S$ admits the positive first Chern class which implies that 
the anti-canonical line bundle on $S$ is positive. 
Then by the same procedure as in K\"ahler manifolds, we obtain the vanishing theorem of Kodaira and the Serre duality of 
the basic cohomology groups. 
The basic cohomology groups inherit many of the same properties as the ordinary Dolbeault cohomology groups
of the K\"ahler manifolds \cite{El K, K-T, BG}. 
Then we have the vanishing of cohomology groups,
\bgn{lemma}\label{lem: Kodaira vanishing}\label{lem: 5.3}
Let $S$ be a compact positive Sasaki manifold of real dimension $2n-1.$ 
Then we have 
$$H^{q,0}_B(S)=H^{0,q}_B(S)=\{0\},\quad \forall q>0.$$
In particular, $H^{2,0}_B(S)=H^{0,2}_B(S)=\{0\}$.
\end{lemma}
\bgn{proof}
Since the transversal complex structure on $S$ admits the positive first basic Chern class which is represented by a positive real basic form of type $(1,1)$, there is a Sasakian structure with the same underlying transversal complex structure whose transversal Ricci curvature is positive
from the transverse Calabi-Yau theorem \cite{El K}, \cite{BGM}. 
Let $\eta$ be the contact form on the $(2n-1)$ dimensional Sasakian manifold $S$. 
Then we have the pairing $\til\Psi $ of the basic forms $\a,\b$ 
by 
$$
\til\Psi(\eta\w\a\w\b):=\int_S\eta\w\a\w\b
$$
Since $d\eta$ is a basic $2$-form, we have the integration by part.
By using the basic Hodge star operator $*_B$, we have the $L^2$-metric
$\|\a\|^2_{L^2}=\int_S\eta\w\a\w *_B\a=\int_S(\a,\a)\vol_S$,
(see Boyer-Galicki,\cite{BG} Sasakian Geometry,  page 217.)
where $(\,,\,)$ is the metric on forms and $\vol_S$ is the volume form on the Sasakian manifold $(S, g)$.
Then the transversal K\"ahler structure on $S$ gives a basic, Hermitian connection $A$ on the transversal holomorphic line bundle $K^{-1}$ such that the curvature form $F_A$ is a positive, basic $(1,1)$ form  on $S$. 
We denote by $\ol\pa_B$ the $\ol\pa$-operator acting on $K^{-1}$-valued basic forms.
Let $\square_B^{\ol\pa}$ be the basic Laplacian of $\ol\pa$-operator acting on $K^{-1}$-valued basic forms
and 
$\square_B^{\pa}$ the complex conjugate of $\square_B^{\ol\pa}$.
Then we have the Kodaira-Nakano identity, 
$$
\square_B^{\ol\pa}=\square_B^{\pa}+[\sqrt{-1}F_A, \w]
$$
(The identity follows from a local calculation as in K\"ahler geometry, see for instance Griffith-Harris, Principle of algebraic geometry, page 154-155.)
Let $\phi$ be a $K^{-1}$-valued basic form of type $(n-1,q)$, $q>0$. Then we have the inequality using the basic Hermitian metric on $K^{-1}$,
$$
\|\ol\pa_B\phi\|^2_{L^2}+\|\ol\pa^*_B\phi\|^2_{L^2}\geq \int_X ( [\sqrt{-1}F_A, \w]\phi, \phi)\vol_S
\geq 2\pi q\|\phi\|^2_{L^2}
$$
Since $K^{-1}$ admits a basic Hermitian connection,
we have the de Rham-Hodge theorem of $K^{-1}$-valued basic forms of type
$(n-1,q)$ with respect to the Dolbeault Laplacian $\square_B^{\ol\pa}$ 
(see \cite{El K}).
Thus every element of cohomology group $H^{n-1,q}_B(S, K^{-1})$ is represented by the $K^{-1}$-valued basic harmonic forms w.r.t $\square_B^{\ol\pa}$
(Note that $S$ is a transversally complex manifold of dimension $n-1$). 
If $\phi$ is harmonic, it follows from the inequality that $\phi=0$.
Then we have the vanishing $H^{n-1,q}_B(S, K^{-1})=\{0\}$. 
There is an isomorphism $H^{n-1, q}_B(S, K^{-1})\cong H^{0,q}_B(S, K\otimes K^{-1})$. 
Since $K\otimes K^{-1}$ is a trivial complex line bundle, we have 
we have the vanishing of the cohomology groups
$$H^{0,q}_B(S)=\{0\}, \qquad \forall q>0.$$
By using the basic Dolbeault decomposition (see also page 217 in \cite{BG}), we have 
$H^{q,0}_B(S) =\ol{H^{0,q}_B(S)}=\{0\},\, (q>0).$
\end{proof}
{\sc Remark.}
Note that the transverse Dolbeault theorem is not necessary to show the vanishing of the cohomology groups of Lemma \ref{lem: 5.3}.
The vanishing of the cohomology groups on Lemma \ref{lem: 5.3} was already shown in \cite{BGM} when the positive Sasaki manifolds are regular or quasi-regular.
This vanishing theorem on positive Sasaki manifolds 
can be proved by the different method (Bochner technique) \cite{Con}.\\
\bgn{lemma}\label{lem: 5.4}
$$H^2_B(S) = H^{1,1}_B(S) = \R dd^c t\oplus H^2(S, \R)$$
\end{lemma}
\bgn{proof}
It is shown that the de Rham cohomology group $H^i(S)$ of the Sasakian manifold $S$ coincides with 
the basic primitive cohomology groups $P_B^i(S)$, for $1\leq i\leq n-1$ (for instance, see the proposition 7.4.13, pp. 233 in \cite{BG}).
It follows from the lemma \ref{lem: Kodaira vanishing} that 
$$
H^2_B(S) =H^{1,1}_B(S) = \R dd^c t\oplus P^2_B(S) = \R dd^c t \oplus H^2(S)
$$
\end{proof}
\bgn{lemma}
Let $\pi : X\to X_0$ be a resolution of an affine variety $X_0$ with only  normal isolated singularity with trivial $K_X$. Then we have the vanishing 
\bgn{equation}\label{eq: GR vanishing}
H^i(X, \O_X) =\{0\}, \qquad  ( \forall \,i>0).
\end{equation}
Further let $C(S)$ be the complement $X\bsh E$ where $E$ is the exceptional set of the resolution $X$. 
Then we also have 
\bgn{equation}
H^i(C(S), \O_{C(S)}) =\{0\} , \qquad ( n-1 > \forall \,i >0).
\end{equation}
In particular, if $n\geq 3$, we have $$ H^1(C(S), \O_{C(S)}) =\{0\},$$ where $\dim_\C X =n$.
\end{lemma}
\bgn{proof}
Since $X_0$ is an affine variety and $H^p(X_0,R^q\pi_*\O_X)=\{0\}$ 
for $p>0$, 
the first vanishing (\ref{eq: GR vanishing}) follows from the Grauert-Riemenschneider vanishing theorem, 
$$
R^i\pi_*K_X=R^i\pi_*\O_X =0, \qquad  ( i>0)
$$
Let $H_E^i(X, \O_X)$ be the local cohomology groups with supports in $E$, and coefficients in the structure sheaf 
$\O_X$. 
Then the short exact sequence,
$$
0\to \Gam_E(X, \O_X) \to \Gam(X, \O_X) \to \Gam (C(S), \O_{C(S)}) \to 0
$$
yields the long exact sequence, 
\bgn{equation}\label{eq: long exact seq}
\cdots \to H^i_E(X, \O_X) \to H^i(X, \O_X) \to H^i( C(S), \O_{C(S)}) \to H^{i+1}_E(X, \O_X) \to \cdots
\end{equation}
By applying the local duality theorem of the cohomology groups to the resolution $\pi: X\to X_0$, we have 
$$
\dim H^j_E(X, \O_X) =\dim R^{n-j}\pi_*K_X =\dim R^{n-j}\pi_*\O_X,
$$ 
(see for instance, corollary 2.5.15,  page 60 \cite{Ishi}).
Then it follows from the Grauert-Riemenschneider vanishing theorem 
 that 
$$H^j_E(X, \O_X) =\{0\}$$ for all $n-j >0$. 
Thus from the exact sequence (\ref{eq: long exact seq}) and (\ref{eq: GR vanishing}), we obtain 
$$
H^i(C(S), \O_{C(S)}) =\{0\}, \qquad ( n-1>\forall\, i >0).
$$
 It implies that $H^1(C(S), \O_{C(S)})=\{0\}$ if $n\geq 3$.
\end{proof}\par
The resolution $\pi: X\to X_0$ gives the identification $\pi |_{X\bsh E} :X\bsh E\cong X_0\bsh\{p\}=C(S)$. 
We denote by  $\{ t>T\}$ the subset of the cone $C(S)=\R\times S$ defined by
$$
\{t>T\}:=\{\, (t, s)\in \R\times S\, |\,  t>T\, \},
$$
 for a constant $T$.
The subset $\{ t>T\}$ is regarded as the subset of $X$  by the identification $\pi |_{X\bsh E} :X\bsh E\cong X_0\bsh\{p\}=C(S)$. 
Then the set $\{ t>T\}$ is an unbounded  region of $X$ and 
we also denote by $\{ t\leq T\}$ 
the complement of $\{ t>T\}$ in $X$ which is a compact set including the exceptional set $E$.
Let $p_{\ss S}$ be the projection from $C(S)=\R\times S$ to $S$. We use the same notation $p_S$ for 
the restricted projection $\{ t>T\} \to S$.
\bgn{lemma}\label{lem: tilkappa}
Let $\k$ be an arbitrary K\"ahler form on $X$ with the K\"ahler class $[\k]\in H^2(X,\R)$. We assume $\dim_\C X\geq 3$.
Then there is a real form $\til \k$ of type $(1,1)$ on $X$ which satisfies the followings, 
\bgn{itemize}
\item[\text{\rm (i)}]  
$[\til \k]=[\k]\in H^2(X, \R)$
\item[\text{\rm (ii)}]
For a sufficiently large constant $T_+$, the restriction $\til\k$ to the subset $\{ t>T_+\}$ is given by the pull back of 
a $d$-closed, primitive basic form $\h\k$ of type $(1,1)$ on $S$, i.e., 
$$
\til\k|_{\{t>T\}} = p_{\ss S}^* \h\k.
$$
\end{itemize}

\end{lemma}
\bgn{proof}
It follows from the lemma \ref{lem: 5.3} and \ref{lem: 5.4} that we have 
$$
H^2(C(S))\cong H^2(S)\cong P^2_B(S)\cong P^{1,1}_B(S)
$$
Then there is a $d$-closed, primitive basic form $\h\k\in P^{1,1}_B(S)$ and a $1$-form $\theta$ on $C(S)$ such that 
\bgn{equation}\label{eq: hkappa}
\k|_{C(S)} =p_S^*\h\k -d\theta.
\end{equation}
Then $d\theta$ is a form of type $(1,1)$ on $C(S)$ where $\theta=\theta^{1,0}+\theta^{0,1}$ and we have 
$$
d\theta=\ol\pa\theta^{1,0}+\pa\theta^{0,1}, \qquad \pa\theta^{1,0}=0,\qquad\ol\pa\theta^{0,1}=0
$$
It follows from the lemma \ref{eq: GR vanishing} that there is a function $\phi$ on $C(S)$ such that
$\theta^{0,1}=\ol\pa \phi$. Since $\theta$ is a real form, it follows from $\theta^{1,0}=\ol{\theta^{0,1}}$ that 
\bgn{align*}
d\theta=&\pa\theta^{0,1}+\ol{\pa\theta^{0,1}}\\
=&\pa\ol\pa\phi+\ol\pa\pa\ol\phi\\
=&2\sqrt{-1}\pa\ol\pa\phi^{Im}
\end{align*}
where $\phi^{Im}=\frac1{2\sqrt{-1}}(\phi-\ol\phi)$ is the imaginary part of $\phi$.
Let $\rho_T$ be a smooth function such that 
$$
\rho_{\ss T}(x) =
\bgn{cases}
&1 \qquad  x\in \{ T<t\}\\
&0\qquad  x\in \{ t \leq T-1\}
\end{cases}
$$
and $0\leq\phi_T(x)\leq 1$ on $X$, i.e., $\rho_T$ is a cut off function which takes the value $0$ on a neighborhood of $E$ and $1$ at the infinity.
We define $\til\k$ by 
$$
\til\k=k + 2\sqrt{-1}\pa\ol\pa \rho_{\ss T_+}\phi^{Im}.
$$
Then we see that $\til\k$ is a form of type $(1,1)$ on $X$ which is the K\"ahler form $\k$ on 
$\{t\leq T_+ -1\}$ and $p_{\ss S}^*\h\k$ on $\{ t>T_+\}$ from (\ref{eq: hkappa}).
Thus $\til\k$ satisfies the conditions.
\end{proof}
\bgn{lemma}\label{lem: initial kahler metric}
Let $\k$ be an arbitrary K\"ahler form on $X$. 
We assume $\dim_\C X \geq 3$. Then there exists a K\"ahler form $\ome_{\k,0}$ which satisfies the followings, 
\bgn{itemize}
\item[\text{\rm(i)}]
$[\ome_{\k,0}]=[\k]\in H^2(X,\R).$
\item[\text{\rm(ii)}]
There are constants $T_+$ and $T_-$ such that 
$\ome_{\k,0}$ coincides with $c\ome_0+\til\k$ on the region $\{1+ T_-< t\}$, 
where $c$ is a positive constant and $\ome_0$ is the K\"ahler cone metric on $C(S)$ as before and 
$\til\k$ is the form in lemma \ref{lem: tilkappa}, where $1+T_-<T_+$
\end{itemize}
\end{lemma}
\bgn{proof}
We see that there is a positive function $\psi(x)$ on $\R$ which satisfies 
$$\psi(x)=
\bgn{cases}
&x\,\qquad\qquad  \text{\rm if }  x>e^{2(1+T_-)}\\
&\text{\rm constant}\quad    \text{\rm if } e^{2T_-} > x
\end{cases}
$$
and $\psi'(x)\geq 0$ and $\psi''(x) \geq 0$  on the region $\{ e^{2T_-}\leq  x\leq e^{2(1+ T_-)}\}$. 
Then the composite function $\psi( r^2)$ is a function on $X$ 
since $\psi$ is a constant on the region $\{ t< T_-\}$, where $r=e^t$ and we have
$$
\sqrt{-1}\pa\ol\pa \psi(r^2)=\sqrt{-1}\psi''(r^2)\pa r^2\w\ol\pa r^2+\sqrt{-1}\psi''(r^2)\pa\ol\pa r^2
$$
It implies that $\sqrt{-1}\pa\ol\pa \psi(r^2)$ is semi-positive on the region $\{ t\leq 1+ T_-\}$ which is the 
K\"ahler cone metric $\ome_0$ on the region $\{ 1+T_- < t\}$. 
Thus we define $\ome_{\k,0}$ by 
$$
\ome_{\k,0}= c \sqrt{-1}\pa\ol\pa \psi(r^2)+ \til\k.
$$
Then $\ome_{\k,0}$ satisfies both conditions (i) and (ii). 
So it suffices to show that $\ome_{\k,0}$ is a K\"ahler form for 
constants $T_+$, $T_-$ and $c$, which are taken by the following. 
We divide $X$ into the following five regions: 

(1) $\{ t<T_-\}\qquad $    (2)$\{ T_-\leq t\leq 1+ T_-\}\qquad$
(3) $\{1+ T_-< t < T_+-1\}$,\par \medskip
(4) $\{T_+-1 \leq t\leq T_+\}\qquad $
(5) $\{T_+ <t\}$\par\noindent
On the region (1), $\ome_{\k,0}$ is the K\"ahler form $\k$ and on the region (3), 
$\ome_{\k,0}$ is $c\ome_0+k$ which is also K\"ahlerian. 
On the region (2), $\ome_{\k,0}= k+\sqrt{-1}\pa\ol\pa \psi(r^2) $ is positive since
$\sqrt{-1}\pa\ol\pa \psi(r^2)$ is semi-positive and $\k$ is positive. 
On the region (5), $\ome_{\k,0}$ is a sum of the K\"ahler cone form $\ome_0$ and  the bounded form $p_{\ss S}^*\h\k$ with respect to the cylinder metric.
Since the cone metric $\ome_0$ grows with order $O(e^{2t})$, there is  a sufficiently large $T_+=T_+(c_0)$  for a positive $c_0$ such that 
$c_0 \sqrt{-1}\pa\ol\pa \psi(r^2)+ \til\k$ becomes positive.
Finally since the region (4) is compact, there is a positive $c$ with $c>c_0$
such that $\ome_{\k,0}$ is a K\"ahler form on the region (4). We can see that $\ome_{\k,0}$ is still positive for $c>c_0$
on (5). 
Hence $\ome_{\k,0}$ is a K\"ahler form.

\end{proof}

\bgn{proof} \negthinspace\negthinspace{\it of theorem} 
\ref{th: existence theorem crepant}
In the case $\dim_{\C}X=2$,  $X$ is a minimal resolution of 
the ordinary double point (see \cite{Kr}).
 Thus every K\"ahler class is represented by 
the class of the exceptional divisors which lies in compactly supported cohomology groups 
and we can have a K\"ahler form $\ome_{\k,0}$ as in the lemma \ref{lem: initial kahler metric}
with $\h k=0$. 
Hence we have the initial K\"ahler form $\ome_{\k,0}$ as in lemma \ref{lem: initial kahler metric} for 
every K\"ahler class on $X$ of dimension $n\geq 2$.
We recall that a positive function $F$ is defined by 
$$\Ome\w\ol\Ome =c_nF_k\ome_{k,0}^n$$ with 
$F_k=e^{f_k}$. 
 We shall show that the function $F_\k$ satisfies 
 $$\|e^{(2+\del)t}(F_k-1)\|_{C^{k,\a}}<\infty,  \qquad \del >0,\quad k\geq 3,\,0<\a<1$$
 with respect to the cylindrical metric  $dt\w d^ct +dd^c t= e^{-2t}\ome_0$ 
 Then the result follows from theorem \ref{th: existence CY on conical}
 
Since $\ome_0$ is a Ricci-flat K\"ahler cone metric, we have 
 $$
 \Ome\w\ol\Ome=c_n\ome_{0}^n
 $$
Then $\ome_0^n=e^{f_k}\ome_{k,0}^n$. 
Since $\ome_{k,0}=\ome_0+\til\k=\ome_0+p_{\ss S}^*\h\k$ on the region $\{t>T_+\}$,
we have 
 \bgn{align}
 e^{-f_k}=&\frac{\ome_{k,0}^n}{\ome_0^n}=\frac{(\ome_0+\til k)^n}{\ome_0^n}\\
 =&1+n\frac{\til k\w \ome_0^{n-1}}{\ome_0^n}
 \label{eq: primitive key point}\\
 +&\sum_{i=2}^n\frac{n!}{i!(n-i)!}\frac{\til{k}^i\w\ome_0^{n-i}}{\ome_0^n}\label{eq: the rest of terms}
 \end{align}
 Since $\til k=p_{\ss S}^*\h\k$ is the pullback of a primitive, basic $(1,1)$-form on $S$ with respect to 
 $dd^ct$, 
 we have 
 $\til k\w(dd^c t)^{n-2}=0$.
 Thus  the second term of  (\ref{eq: primitive key point}) vanishes since we have
 \bgn{align}
 (n-1)^{-1}\til k\w\ome_0^{n-1}=&\til k\w e^{2(n-1)t} (2dt\w d^ct\w(dd^c t)^{n-2})\\
 =&e^{2(n-1)t} (2dt\w d^ct\w \h\k\w(dd^c t)^{n-2})=0
 \end{align}
 
Each term of (\ref{eq: the rest of terms}) is given by 
\bgn{align}
(n-i)^{-1}\frac{\til{k}^i\w\ome_0^{n-i}}{\ome_0^n}=&
\frac{\til k^i\w e^{2(n-i)t} (2dt\w d^ct\w(dd^c t)^{n-i-1})}{
e^{2nt} (2dt\w d^ct\w(dd^c t)^{n})}\\
=&e^{-2i t}\frac{\til k^i\w(2dt\w d^ct\w(dd^c t)^{n-i-1})}
{(2dt\w d^ct\w(dd^c t)^{n})}
\end{align}
Since the $1$-form $dt$ and $d^ct$ and the $2$-form $dd^c t$, $\til k$ are bounded with respect to the cylinder metric, we have 
$$
\frac{\til{k}^i\w\ome_0^{n-i}}{\ome_0^n}=O(e^{-2it})=O(r^{-2i})
$$
Note that $\til\k=p^*_S\h\k$ for a basic $2$-form $\h\k$.
It follows from $i\geq 2$ that 
the every term in (\ref{eq: the rest of terms}) decays with the order 
$O(e^{-4t})$. Hence we have $e^{-f_\k}-1=O(e^{-4t})$. 
Since $e^{-f_\k}$ is bounded, 
we obtain 
$$e^{f_k}-1 =O(e^{-4t}).$$

We also have the estimate of the higher order derivative of $F_k$ since 
we use the $C^k$-norm with respect to the cylinder metric.
Note that the $C^k$-norm of $e^{-4t}$ is just estimated the derivative by 
the cylinder parameter $t$ which also decays with the order $O(e^{-4 t})$.

Hence there is a solution $u$ of the \M-A equation
$$\Ome\w\ol\Ome=c_n(\ome_{k,0}+dd^cu)^n$$ from the existence theorem 
\ref{th: existence CY on conical}.
Then the K\"ahler form $\ome_k=\ome_{k,0}+dd^cu$
gives a Ricci-flat K\"ahler metric on $X$
\end{proof}\par
\bgn{remark}
If  $\ome_{\k,0}$ is not in the form $\ome_0+p_S^*\h\k$ for the pullback of the basic, primitive $(1,1)$ form $\h\k$
in the region $\{t>T_+\}$, 
$F-1$ only decays with the order $O(e^{-2 t})$ which we can not apply 
the existence theorem \ref{th: existence CY on conical}
\end{remark}

\section{Examples of Calabi-Yau structures on crepant resolutions}
We construct examples of Ricci-flat conical K\"ahler metrics.
Some of them are already known, however there are new  
Ricci-flat conical K\"ahler metrics included
whose K\"ahler classes do not belong to  the compactly supported cohomology groups. 
We start with the trivial example, 
\bgn{example}
The complement $\C^n\bsh \{0\}$ is the cone of the sphere of dimension 
$2n-1$ and 
the standard K\"ahler metric $\ome_{st}$ on $\C^n$ is a conical Ricci-flat K\"ahler metric  
with $r^2=\sum_{i=1}^n|z_i|^2$. 
The induced metric on the sphere is an Einstein-Sasakian metric.
\end{example}
\bgn{example}
Let $\Gam$ be a finite subgroup of the special unitary group SU$(n)$ which freely acts on 
$\C^n\bsh\{0\}$. Then the quotient $X_0=\C^n/\Gam$ has a normal isolated singularity at the origin $0$ and 
the complement $\C^n\bsh\{0\}/\Gam $ is the cone $C(S)$ of the 
Einstein-Sasakian manifold $S:=S^{2n-1}/\Gam$. 
The induced metric on the cone $C(S)$ is a Ricci-flat K\"ahler cone metric. 
Then the theorem \ref{th: existence theorem crepant} shows that a crepant resolution of the isolated quotient 
singularity $\C^n/\Gam$ has a Ricci-flat conical K\"ahler metric. 
This class is already obtained by Kronheimer \text{\rm\cite{Kr}} for $n=2$ and by Joyce \text{\rm \cite{Jo}} for $n>2$.
\end{example}
\bgn{example}
Let $Z$ be a compact  K\"ahler-Einstein manifold with positive 1-st Chern class $c_1(Z)$, which is called a Fano manifold 
(We assume that $\dim_{\C} Z=n-1\geq 2$).
Then the sphere bundle $S$ of the canonical line bundle $K_Z$ of the Fano manifold $Z$ admits an Einstein-Sasakian metric and the complement 
of the zero-section $K_Z\bsh\{0\}$ is the cone $C(S)$ which has the 
Ricci-flat K\"ahler cone metric. 
Then we apply the theorem \ref{th: existence theorem crepant} to $X=K_Z$ and obtain a Ricci-flat conical K\"ahler metric in every k\"ahler class of $K_Z$. 
Calabi \text{\rm \cite{Cal}} already constructed a Ricci flat K\"ahler metric on $K_Z$ by the bundle construction whose K\"ahler class lies in the compactly supported cohomology group, i.e., the anti-canonical class. 
Since we have the vanishing of the cohomology groups 
$H^1(S)=\{0\}$ and $H^3_{\cpt}(X)\cong H^3(X, X\bsh Z) \cong H^1(Z)=\{0\}$ 
by the duality theorem. 
Then we obtain the exact sequence, 
$$
0\to H^2_{cpt}(X)\to H^2(X)\to H^2(S)\to 0.
$$
We also have $H^2_{\cpt}(X)\cong H^2(X, X\bsh Z) \cong H^0(Z)\cong \R$.
Thus $\dim H^2(X)=\dim H^2(S)+1$.
If $b_2(Z)=\dim H^2(Z)=\dim H^2(X)$ is greater than $1$, it follows from the exact sequence that there is a K\"ahler class which does not belong to the 
compactly supported cohomology group. 
Thus in the cases, we obtain a new family of Ricci-flat K\"ahler metrics on $K_Z$. 
\end{example}
\bgn{example}\label{ex: toric Sasakian}
Recently Futaki-Ono-Wang \text{\rm\cite{FOW}} obtained an Einstein-Sasakian metric on the sphere bundle of the canonical line bundle of every toric Fano manifold which gives the Ricci-flat K\"ahler cone metric on the complement $K_Z\bsh\{0\}$. 
Thus the theorem 5.1 shows that there exists  a Ricci-flat conical K\"ahler metric in every K\"ahler class on the total space of the canonical line bundle on 
an arbitrary toric Fano manifold.
If the canonical line bundle is the $m$-th tensor of a line bundle $L$ on a toric Fano manifold,
then we can apply the theorem to the total space $L$ to obtain 
a Ricci-flat conical K\"ahler metric on $L$ in each K\"ahler class.
Let $\widehat{\ CP^2}$ be the blown up $\C P^2$ at one point.
Then $\widehat{\ CP^2}$ is a toric Fano manifold which does not admit a K\"ahler-Einstein metric. 
However it is shown that the canonical line bundle of 
$\widehat{\ CP^2}$ admits a family of Ricci-flat conical K\"ahler metrics 
which is parametrized by the open set $H^2( \widehat{\ CP^2})$, i.e., the K\"ahler cone of $Z$.
Note that  $\dim H^2(\widehat{\ CP^2})=2$
and $\dim H^2_{cpt}(\widehat{\ CP^2}) =1.$
Sasaki-Einstein metrics on the blown up $\C P^2$ at one point 
was explicitly constructed by Martelli and Sparks \cite{MS1}
and Oota-Yasui \cite{OY} described a Ricci-flat metric on 
a resolved Calabi-Yau cone. 
\end{example}
\bgn{example}
We denote by $\Sigma_k$ a blown up $\C P^2$ at $k$ points in general position \text{\rm($0\leq k \leq 8$)}. 
Then it is known that $\Sigma_k$ admits an Einstein-K\"ahler metric for $3\leq k\leq 8$ \cite{TY3}, \cite{Ti1}.
Since $\Sigma_k$ for $k=1,2$ is a toric Fano surface, there is a Einstein-Sasakian metric on the sphere bundle of the canonical line bundle 
on $\Sigma_k$ \cite{FOW}, \cite{MSY}.
Hence we obtain a complete Ricci-flat K\"ahler metrics in every K\"ahler class of the canonical line bundle 
of a blown up $\C P^2$ at generic $k$ points for all $0\leq k \leq 8$.
\end{example}
\bgn{example}
Let $X_0$ be the $3$ dimensional hypersurface in $\C^4$ defined by 
$z_0^2+z_1^2+z_2^2+z_3^2=0$. 
The hypersurface $X_0$ has a singularity at the origin $0$ which is the ordinary double point and 
there is a small resolution $X\to X_0$ with trivial $K_X$. 
The small resolution $X$ is the total space $\O(-1)\oplus \O(-1)$ on $\C P^1$ and $X_0$ is obtained by the contraction of the zero-section to one point.
The complement $X_0\bsh\{0\}$ is the cone $C(S)$ of the sphere bundle $S$
over the product $\C P^1\times \C P^1$. 
Thus the sphere bundle $S$ has an Einstein-Sasakian metric which induces the 
Ricci-flat K\"ahler cone metric on $C(S)$. 
Then we obtain a one dimensional family of Ricci-flat conical K\"ahler metrics on the resolution,
since $\dim H^2(X)=1$ and $H^2_{\cpt}(X)=0$. 
Candelas and de la Ossa \text{\rm \cite{Can}} described  a Ricci-flat metric on 
a conifold
and it is also keen to show that the Ricci-flat metric in \text{\rm \cite{Can}}
coincides with the one in our family constructed by the theorem 5.1.
\end{example}

\section{Appendix}
 We shall show the following lemmas \ref{lem: from weighted CY to Einstein -Sasakian } and \ref{lem: from Einstein-Sasakian to weighted CY} to  prove the proposition \ref{prop: Einstein-Sasakian and weighted Calabi-Yau}. 
 Let Ric$_{g_{\ss S}}$ be the Ricci curvature of a Sasakian manifold 
 $(S, g_{\ss S})$. 
We see that  $d\eta_{\ss S}$ is a K\"ahler form on the distribution $D$ which satisfies
 $\L_{\xi_{\ss S}} d\eta_{\ss s}=0$.
 Then $d\eta_{\ss S}$ gives the transversal Riemannian metric $g^T_{\ss S}$ on $D$ with 
 $\L_{\xi_{\ss S}} g^T_{\ss S}=0$.
 Let Ric$_{g^T_{\ss S}}$ be the Ricci curvature of the transversal Riemannian metric $g^T_{\ss S}$. 
 Then there is a relations between 
  Ric$_{g^T_{\ss S}}$ and  Ric$_g$, 
 \bgn{proposition}\label{prop: transversal Ricci} (Boyer-Galicki, pp 224, theorem 7.3.12)
 \bgn{align*}
 \text{\rm Ric}_g (u,v)= &\text{\rm Ric}_{g^T_{\ss S}}(u,v)-2g(u,v), \qquad u, v\in D \\
 \Ric_g(u,\xi) =& 2(n-1)\eta_{\ss S}(u),\qquad\forall \,u \in TS
 \end{align*}
 \end{proposition}
 Thus a Sasakian metric $g_S$ is Einstein if and only if the transversal metric 
 $g^T$ is Einstein with Einstein constant $2n$.
 
 \bgn{lemma}\label{lem: from weighted CY to Einstein -Sasakian }
 Let $(\Ome, \ome)$ be a weighted Calabi-Yau structure on $C(S)$ and 
 $(\eta_{\ss S}, \psi_{\ss S})$ the corresponding structure defined in \text{\rm (\ref{ali: psi, eta})} which induces a Sasakian structure on $S$. 
 Let  $g_{\ss S}^T$ be the transversal Riemannian metric given by 
 the transversal K\"ahler structure  $d\eta_{\ss S}$ on $S$.
 Then $g_{\ss S}^T$ is Einstein with Einstein constant $2n$, that is, Ric$_{g^T_{\ss S}}= 2n g_{\ss S}^T$.
 \end{lemma}
 \bgn{proof}
 We take an open covering $\{U_\a\}_\a$ of $S$ such that 
 each $U_\a$ admits coordinates $(x, z_1^\a,\cdots, z_{n-1}^\a)$ which is compatible with transversal holomorphic structure on $S$, that is, 
$(z_1^\a,\cdots, z_{n-1}^\a)$ gives the transversal complex structures and 
 $\frac{\pa }{\pa x} =\xi_{\ss S}$.
 (Such coordinates are called foliation coordinates.)
 Then we have a $d$-closed, local holomorphic $n$-form $\Ome_{S,\a}=dz^\a_1\w\cdots \w dz_{n-1}^\a$ which is a basic section of $K_D|_{U_\a}$.  
 For simplicity, we write $\Ome_{S,\a}$ for its pullback 
 $p_{\ss S}^*\Ome_{S,\a}$ on $C(S)$.
 Since $\psi= i_{\frac{\pa }{\pa t}}\Ome$ is a holomorphic section of  $p_{\ss S}^*K_D$, there is a holomorphic function $f_\a$ on $p_{\ss S}^{-1}(U_\a)$ such that 
 $$
 \psi=e^{f_\a}\Ome_{S,\a}.
 $$
Substituting $ \psi=e^{f_\a}\Ome_{S,\a}$ into (\ref{ali: r-Sasaki-M-A}),
 we have 
 \bgn{align}\label{ali: p a }
 e^{f_\a+\ol f_\a}\eta\w\Ome_{S,\a}\w\sig(\ol\Ome_{S,\a})=c_{n-1}r^{2n}\eta\w(\frac12d\eta)^{n-1}
 \end{align}
 
 We define a function $k_\a$ locally by
 \bgn{equation}\label{eq: k a }
 \Ome_{S,\a}\w\sig(\ol\Ome_{S,\a})=e^{k_\a}
 c_{n-1}(\frac12d\eta_{\ss S})^{n-1}
\end{equation}
Then the Ricci form of the transversal K\"ahler form $g_{\ss S}^T$ on an open set $U_\a$ of $S$ is given by 
Ric$_{g_{\ss S}^T}=\sqrt{-1}\pa_B\ol\pa_Bk_\a$, where 
$\olpa_B$ is the 
$\ol\pa$-operator with respect to the transversal complex structure on $S$ and $\pa_B$ is its complex conjugate.
By multiplying $\eta=\eta_{\ss S}$ to both sides of (\ref{eq: k a }) and pulling back by $p_{\ss S}$ to $C(S)$, we have 
\bgn{equation}\label{eq: k a -2}
 \eta\w \Ome_{S,\a}\w\sig(\ol\Ome_{S,\a})=e^{p_{\ss S}^*k_\a}c_{n-1}\eta\w(\frac12d\eta_{\ss S})^{n-1}
\end{equation}
Comparing to (\ref{ali: p a }), we have 
\bgn{equation}\label{eq: 4-24}
p_{\ss S}^*e^{k_\a}= r^{2n}e^{-(f_\a+\ol f_\a)}=e^{2nt}e^{-(f_\a+\ol f_\a)}
\end{equation} 
A function $p_{\ss S}^*\k_\a$ on  $C(S)$ is basic with respect to both $\xi, \frac{\pa}{\pa t}$.
Thus 
we have $p_{\ss S}^*\pa_B\ol\pa_B k_\a=\pa\ol\pa p_{\ss S}^*\k_\a$.
The it follows that
$$
\text{\rm Ric}_{g_{\ss S}^T} =\sqrt{-1}\pa_B\ol\pa_B k_\a=i_{\ss S}^*\sqrt{-1}\pa\ol\pa p_{\ss S}^*\k_\a,
$$
where $\ol\pa $ is the $\ol\pa$-operator with respect to the complex structure on $C(S)$.
From (\ref{eq: 4-24}), we have 
\bgn{align}
\text{\rm Ric}_{g_{\ss S}^T}=&2ni_{\ss S}^*\sqrt{-1}\pa\ol\pa t -i_{\ss S}^*\sqrt{-1}\pa\ol\pa ( f_\a+\ol f_\a)
\end{align}
Since $\ol\pa f_\a=0$, we have $\pa\ol\pa ( f_\a+\ol f_\a)=0$ and
the transversal K\"ahler form is  $\frac12d\eta_{\ss S}=\sqrt{-1}\pa\ol\pa t|_{\ss S}$.
Thus we obtain 
\bgn{align}
\text{\rm Ric}_{g_{\ss S}^T}=&2n\sqrt{-1}\pa\ol\pa t =
2n\, (\frac12d\eta_{\ss S})
\end{align}
Then it follows that  Ric$_{g_{\ss S}^T}= 2n\, g_{\ss S}^T$.
\end{proof}\par
Conversely the following lemma shows that 
an Einstein-Sasakian structure on $S$ gives the weighted Calabi-Yau structure on $C(S)$, 
\bgn{lemma}\label{lem: from Einstein-Sasakian to weighted CY}
An Einstein-Sasakian structure on $S$ corresponds to the weighted Calabi-Yau structure on $C(S)$ under the correspondence in the proposition \ref{prop: Sasakian and kahler cone metric}
\end{lemma}
\bgn{proof}
As in the proposition \ref{prop: transversal Ricci}, 
an Einstein-Sasaki structure gives the transversal Einstein-K\"ahler metric with scalar curvature $2n$,
Ric$_{g_{\ss S}^T}= 2n g_{\ss S}^T$.
Thus the $1$-st Chern class $c^1_B(S)$ of the transversal canonical line bundle $K_D$
is represented by the form 
$2n( \frac12d\eta_{\ss S})=n\,d\eta_{\ss S}$.
There is an exact sequence on $S$,
$$
\to H^0(S) \overset{i}{\to} H^2_B(S) \overset{j}{\to} H^2(S)\to 
$$
in which  we have $i(\a)=\a d\eta_{\ss S}$ ( $\a\in H^0(S)$) and 
$j(c^1_B(S)) = n j ( d\eta_{\ss S}) = i\circ j(n)=0$.
Thus $c^1(S)=c^1(K_D)$  vanishes.
Since $S$ is simply connected, the line bundle $K_D$ is trivial.
Let 
$\{U_\a\}$ be an open covering of $S$ such that each $U_\a$ admits 
a foliation coordinates
$(x^\a,z_1^\a,\cdots, z_{n-1}^\a)$.
We denote by  $p_{\ss S}: 
C(S) \to S$  the projection. 
We take $U_\a$  sufficiently small such that 
every $d$-closed holomorphic $1$-form on $p_{\ss S}^{-1}(U_\a)\cong \R\times U_\a$ is written as a $d$-exact $1$-form of a holomorphic function, 
that is, the holomorphic de Rham theorem holds on $U_\a$.
Taking a $d$-closed, basic, holomorphic $(n-1)$-form $\Ome_{S,\a}=dz^\a_1\w\cdots\w dz_{n-1}^\a$ on  $U_\a$, we define a basic function $k_\a$ by  
 
\bgn{align}\label{ali: k a}
\Ome_{S,\a}\w\sig(\ol\Ome_{S,\a})=c_{n-1}e^{\k_\a}(d\eta_{\ss S})^{n-1}
\end{align}
Then the transversal Ricci form Ric$_{g_S^T}$ 
is given by 
\bgn{align}
\text{\rm Ric}_{g_S^T} =\sqrt{-1}\pa_B\ol\pa_B\k_\a
\end{align}
 Since the transversal K\"ahler form  $\frac12d\eta_{\ss S}$ is 
 Einstein-K\"ahler, we have 
\bgn{align}
\sqrt{-1}\pa_B\ol\pa_B\k_\a=&2n(\frac 12d\eta_{\ss S})=
2ni_{\ss S}^*(\sqrt{-1}\pa\ol\pa t)
\end{align}
Pulling back them by the projection $p_s : C(S)\to S$ to $C(S)$,
since $\k_\a$ is a basic function, we have 
\bgn{align}
p_{\ss S}^*\text{\rm Ric}_T(\frac12d\eta_{\ss S}) =&p_{\ss S}^*\sqrt{-1}\pa_B\ol\pa_B\k_\a\\
=&\sqrt{-1}\pa\ol\pa p_{\ss S}^*\k_\a =2n(\sqrt{-1}\pa\ol\pa t)
\end{align}
Then on $p_{\ss S}^{-1}(U_\a)\subset C(S)$, we have 
$$
\sqrt{-1}\pa\ol\pa(p_{\ss S}^*\k_\a-2n t)=0
$$
Since  $U_\a$ is sufficiently small such that every $d$-closed 
holomorphic $1$-form is a $d$ exact form of a holomorphic function, 
there is a holomorphic function $q_\a$ on
$p^{-1}(U_\a)$ satisfying 
$$
\k_\a-2nt =q_\a+\ol q_\a
$$
Then from (\ref{ali: k a}), we have 
\bgn{align}
(e^{-q_\a}\Ome_{S,\a})\w\sig(e^{-\ol q_\a}\ol \Ome_{S,\a}) 
=c_{n-1}r^{2n} (\frac12 d\eta_{\ss S})^{n-1}
\end{align}
Note that $r= e^t$.
We define $\psi$ by $e^{-q_\a}\Ome_{S,\a}$. Then the pair ($\eta_{\ss S}, \psi$) 
satisfies (\ref{ali: r-Sasaki-M-A on S}).
If we set 
$\Ome_\a=(dt+\sqrt{-1}Jdt)\w (e^{-q_\a}\Ome_{S,\a})$,
then $\Ome_\a$ is a local holomorphic $n$-form on $C(S)$ which satisfies
 $$
 \Ome_\a\w\sig(\ol\Ome_\a)= c_n\ome^n
 $$
 where $c_n=\frac{(2\sqrt{-1})^n}{n!}$. 
 Hence the Ricci curvature of the K\"ahler metric $\ome$ on 
 $C(S)$ vanishes.
 Further since $S$ is simply connected, the canonical line bundle on $C(S)$ is trivial. Then there is a holomorphic $n$-form $\Ome $ on $C(S)$ 
satisfying the \M-A equation, 
$$
 \Ome\w\sig(\ol\Ome)= c_n\ome^n
 $$
As we see in the section $1$,

$$
\L_{\frac{\pa}{\pa t}}J=0,\qquad \L_{\frac{\pa}{\pa t}}\ome =2\ome
$$

Since 
${\frac{\pa}{\pa t}}$ is the real part of a holomorphic vector field on $C(S)$, 
there is a holomorphic function  $h$ such that
$$
\L_{\frac{\pa}{\pa t}}\Ome=h\Ome
$$
By the action of the Lie derivative 
 $\L_{\frac{\pa}{\pa t}}$ on both sides of the \M-A equation, we have 
$$
(k+\ol k)\Ome\w\sig(\ol\Ome) =2n c_n\ome^n
$$
It follows that $k+\ol k=2n$. Since $k$ is holomorphic, it follows that $k$ is 
a constant $n$.
Hence 
$$
\L_{\frac{\pa}{\pa t} }\Ome =n \Ome 
$$
Thus $(\Ome, \ome)$ is a weighted Calabi-Yau structure on $C(S)$.
Hence we have the result.
\end{proof}\par
\bgn{proof}{\it of theorem \ref{prop: Einstein-Sasakian and weighted Calabi-Yau}}
The theorem follows from the lemmas \ref{lem: from weighted CY to Einstein -Sasakian } and \ref{lem: from Einstein-Sasakian to weighted CY}
and the proposition \ref{prop: transversal Ricci}.
\end{proof}

\bgn{thebibliography}{99}
\bibitem{Au}
T.~Aubin,
{\it Nonlinear Analysis on Manifolds. \M-A equation},
Grunlehren der mathematischen wissenschaften 252, 
a Series of Comprehensive Studies on Mathematics, 
Springer-Verlag, 1982
\bibitem{BK1} 
S.~Bando and R.~Kobayashi,
{\it Ricci-flat K\"ahler metrics on affine algebraic manifolds}, 
In T.~Sunada, editor Geometry and Analysis on Manifolds, volume 1339 of 
Lecture Notes in Mathematics, pages 20-31, springer Verlag, 1988

\bibitem{BK2}
S.~Bando and R.~Kobayashi,
{\it Ricci-flat K\"ahler metrics on affine algebraic manifolds II}, 
Mathematische Annalen, 287, 175-180, 1990
\bibitem {BG}
C.~Boyer and K.~Galicki,
{\it Sasakian Geometry},
Oxford Mathematical Monographs, Oxford University Press, 
Oxford, 2008 
\bibitem{BGM}
C.~Boyer, K.~Galicki and M.~ Nakamaye,
{\it On positive Sasakian geometry},
Geometriae Dedicata 101: 93-102, 2003
\bibitem{Cal}
E. Calabi, 
{\it M\'etrques k\"ahl\'eriennes et fibr\'es holomorphes}, 
Ann.Sci.\'Ecole Norm. Sup. (4), 12 (2), 269-294, 1979

\bibitem{Can}
P.~Candelas and Xenia~C.~de la Ossa, 
{\it Comments on conifolds},
Nuclear Physics. B, 343(1), 246-268, 1979
\bibitem {Coe}
C.~van Coevering,
{\it Ricci-flat K\"ahler metrics on crepant resolutions of K\"ahler cones}
arXiv: 0806.3728
\bibitem{Coe2}
C.~van Coevering,
{\it Regularity of asymptotically conical Ricci-flat K\"ahler metrics},
arXiv:0912.3946
\bibitem{Con}
Ronan J.~Conlon, 
{\it A vanishing theorem for positive Sasaki manifolds}, 2011
preprint 

\bibitem{El K}
El Kacimi-Alaoui, Aziz,{\it  Op\'erateurs transversalement elliptiques sur un feuilletage riemannien et applications}, [Transversely elliptic operators on a Riemannian foliation, and applications], Compositio Math. 73 (1990), no. 1, 57--106, 58G05 (57R30)
\bibitem{Ev}
L.C.~Evans,
{\it Classical solutions of fully nonlinear, convex, second order elliptic equations},
Comm. Pure Appl. Math., 25, 1982, 333-363
\bibitem{FOW}
A.~ Futaki, H.~Ono and G.~ Wang,
{\it Transverse K\"ahler geometry of Sasaki manifolds and toric Sasaki-Einstein manifolds},
arXiv:math/0607586

\bibitem{Gal}
S.~Gallot, 
{\it A Sobolev inequality and some geometric applications},
Spectra of Riemannian manifolds, Kaigai Publications, Tokyo, 1983, 45-55
\bibitem{Go1}
R.~Goto, 
{\it On hyper-K\"ahler manifolds of type $A_\infty$ and $D_\infty$}, 
Comm. Math. Physics, 
198, No. 2, 469-491, (1998)
\bibitem{G-T}
D.~Gilberg and N.S.~Trudinger, 
{\it Elliptic Partial Differential Equations of Second Order}
Grundlehren der mathematishen Wissenshchaften 224, 
ASeries of Comprehensive Studies in Mathematics, 
Springer-Verlag, 1983
\bibitem{HKLR}
N.J.~Hitchin, A.~Karlhede, U.~Lindstr\"om and M.~R\v{o}cek,
{\it HyperK\"ahler metrics and Supersymmetry, }
Comm. Math. Physics, 108, 535-589, 1987
\bibitem{Ishi}
S.~Ishii, 
{\it Introduction to singularities (Tokuiten nyumon)}, 
in Japanese, 
Springer Verlag, Tokyo, 1997
\bibitem{Jo}
D.~Joyce,
{\it Compact Manifolds with Special Holonomy},
Oxford Mathematical Monograph, Oxford University Press, Oxford, 2000

\bibitem{Kr}
P.~Kronheimer,
{\it The construction of ALE spaces as hyperK\"ahler quotients},
Journal of Diff. Geometry 29, 665-683, 1989
\bibitem{Kry}
N.V.~Krylov, 
{\it Boundedly inhomogeneous elliptic and parabolic equations}, 
Izvestia Akad. Nauk SSSR. 46, 1982, 487-523 [Russina], 
English translation in Math. USSR. Izv., 20, 1983, 459-492
\bibitem{K-T}
F.~Kamber W, P.~Tondeur, {\it de Rham Hodge theory for Riemannian foliations}, Math. Ann. 277 (1987), no. 3, 415--431 53C12
\bibitem{LMc}
R.B.~Lockhart an R.C.~McOwen, 
{\it Elliptic differential operators on noncompact manifolds},
Ann.Scupla Norm.Sup.Pisa Cl. Sci. 12, 1985, pp409-447

\bibitem{Mel}
R.B.~Melrose, 
{\it Atiyah-Patodi-
Singer Index Theorem}, 
A K Peters Ltd, Wellesley, MA, 1983

\bibitem{MS}
D.~Martelli and J.~Sparks, 
{\it Symmetry-breaking vacua and baryon condensates in ADS/CFT}, 
arXiv: physics. 0804.3999, Phys. Rev D79: 065009, 2009
\bibitem{MS1}
D.~Martelli and J.~Sparks, 
{\it Toric geometry, Sasaki-Einstein manifolds and a new infinite class of AdS/CFT duals}, Comm. Math. Phys. 262 (2006), no. 1, 51--89.
\bibitem{MSY}
D.~Martelli, J.~ Sparks and S.~T.~Yau,
{\it  Sasaki-Einstein manifolds and volume minimisation}.
Comm. Math. Phys. 280 (2008), no. 3, 611--673.
\bibitem{Na}
H.~Nakajima, 
{\it Nonlinear Analysis and Complex Geometry},
Iwanami Shoten, Gendai suugaku no tenkai,
(in Japanese)
1999

\bibitem{OY}
T.~Oota and Y.~Yasui, 
{\it Explicit Toric Metric on Resolved Calabi-Yau Cone},
hep-th/0605129,
Phys.Lett. B639, 54-56, 2006
\bibitem{San}
B.~Santoro,
{\it Existence of complete Ricci-flat metrics on crepant resolutions},
math.DG/0902.0595
\bibitem{Siu}
Y.T.~Siu
{\it Lectures on Hermitian-Einstein metrics for stable bundles and K\"ahler-Einstein Metrics}
DMV Seminar, Band 8, Birkh\"auser, 1987

\bibitem{Ti}
G.~Tian, 
{\it Canonical Metrics in K\"ahler Geometry}, 
Lectures in Mathematics, ETH  Z\"urich, 
Birkh\"auser,2000
\bibitem{Ti1}
G.~Tian, 
{\it On Calabi's conjecture for complex surfaces with positive first Chern class}, Invent. Math. 101 (1990), no. 1, 101--172
\bibitem{TY}
G.~Tian and S.-T.~Yau, 
{\it Existence of K\"ahler-Einstein metrics on complete K\"ahler manifolds and their applications to algebraic geometry},
In S.-T.~Yau, editor, Mathematical Aspects of String Theory, volume 1 of advanced series in Mathematical Physics, pages 574-628, World Scientific,
1987
\bibitem{TY1}
G.~Tian and S.-T.~Yau, 
{\it Complete K\"ahler manifolds with zero Ricci curvature. I}
Journal of the American Mathematical Society, 3, 579-609,1980
\bibitem{TY2}
G.~Tian and S.-T.~Yau, 
{\it Complete K\"ahler manifolds with zero Ricci curvature. II}
Inventiones mathematicae, 61, 251-265, 1980
\bibitem{TY3}
G.~Tian and S.-T.~Yau, 
{\it K\"ahler-Einstein metrics on complex surfaces with $C_1>0$}, Comm. Math. Phys. 112 (1987), no. 1, 175--203
\bibitem{Tr}
N.S.~ Trudinger,
{\it Lectures on nonlinear elliptic equations of second order}, 
Lectures on Mathematical Sciences, the University of Tokyo, vol 9, 1995
 
\end{thebibliography}


\address{
Deaprtment of Mathematics \\ 
Graduate School of Science \\
Osaka University\\
Toyonaka, Osaka, 560\\
Japan
}
{goto@math.sci.osaka-u.ac.jp}

\end{document}